\documentclass[english,12pt]{amsart}
\usepackage{amsmath,amsthm,amssymb, amstext, amsopn, amsxtra}
\newtheorem{theorem}{Theorem}[section]
\newtheorem{lemma}[theorem]{Lemma}
\newtheorem{corollary}[theorem]{Corollary}

\newtheorem{prop}[theorem]{Proposition}

\theoremstyle{definition}
\newtheorem{definition}[theorem]{Definition}
\theoremstyle{definitions}
\newtheorem{definitions}[theorem]{Definitions}
\newtheorem{example}[theorem]{Example}
\newtheorem{examples}[theorem]{Examples}

\theoremstyle{remark}
\newtheorem{remark}[theorem]{Remark}
\theoremstyle{remarks}
\newtheorem{remarks}[theorem]{Remarks}
\numberwithin{equation}{section}
\newcommand{\bs}{\bigskip}
\newcommand{\ms}{\medskip}

\newcommand{\sk}{\smallskip}
\newcommand{\D}{\Delta}

\newcommand{\De}{\Delta}
\providecommand\rk{\mathrm{rk\,}}
\providecommand\dim{\mathrm{dim\,}}
\providecommand\im{\mathrm{Im\,}}

\providecommand\deg{\mathrm{deg\,}}
\providecommand\diag{\mathrm{diag\,}}
\providecommand\id{\mathrm{id}}
\providecommand\lspec{\mathrm{lspec}}
\providecommand\rspec{\mathrm{rspec}}
\providecommand\Spec{\mathrm{Spec}}

\begin{document}

\centerline{{\bf {\Large Wedderburn Polynomials over } } }%

\vspace{3mm}

\centerline{{\bf {\Large Division Rings, II} } }

\vspace{7mm}

\centerline{ T. Y. LAM, A.LEROY and A. OZTURK }

\vspace{7mm}

\par
{\small Abstract: A  polynomial $f(t)$ in an Ore extension
$K[t;S,D]$ over a division ring $K$ is a  Wedderburn polynomial if
$f(t)$ is monic and is the minimal polynomial of an algebraic
subset of $K$. These polynomials have been studied in \cite{LL5}.
In this paper, we continue this study and give some applications
to triangulation, diagonalization and eigenvalues of matrices over
a division ring in the general setting of $(S,D)$-pseudo-linear
transformations. In the last section we introduce and study the
notion of $G$-algebraic sets which, in particular, permits
generalization of Wedderburn's theorem relative to factorization
of central polynomials.}

%\footnote{ 1991 {\it Mathematics subject classification}  Primary:
%16D40, 16E20, 16L30; Secondary: 16D70, 16E10, 16G30}

\vspace{10mm}

\section {\bf Introduction}

This paper continues the study of Wedderburn polynomials started
in \cite{LL5}.  Wedderburn polynomials are least left common
multiple of linear polynomials of the form $t-a$ in (skew)
polynomial rings over division rings.  They can be factorized
linearly using Wedderburn's method and have been intensively
studied recently (Cf.\cite{DL},
\cite{LL4},\cite{LL5},\cite{Ro1},\cite{Ro2},\cite{RS1},
\cite{RS2},\cite{Se}).  They appear sometimes under other names
such as rings with separate zeros or polynomials with zeros in
generic positions (\cite{Tr},\cite{GGRW},\cite{GR},\cite{GRW}).
Wedderburn polynomials are also special instances of more general
polynomials called fully reducible (Cf. \cite{Co2},\cite{LO}).

Let us now briefly describe the content of the paper.  In the
sequel $R$ stands for an Ore extension $R=K[t;S,D]$ where $K$ is a
division ring, $S$ an endomorphism of $K$ and $D$ is a
$S$-derivation of $K$.  In section 2 we recall some basic facts
and notations from our previous paper (\cite{LL5}).   In the third
section we present various relations involving the rank of
algebraic sets and, using these, we recover some of the features
of Wedderburn polynomials presented in our previous work. Section
four is devoted to companion matrices. They show up naturally in
the study of the action of $t.$ on $R/Rf$ and are very useful tool
while we characterize when a product of $W$-polynomials is again a
$W$-polynomial.  This generalizes the $(S,D)$-metro equation from
\cite{LL5}.  In section 5 we analyse the problems of
diagonalization and triangulation of matrices over a division
ring.  We work in the general $(K,S,D)$-setting as described
above.  We first study the case of a companion matrix and then,
supposing $S\in Aut(K)$, we analyse the case of a general square
matrix via the companion matrices of its invariant factors.  In
particular we will show that a square matrix $A\in M_n(K)$ is
$(S,D)$-diagonalizable (resp. $(S,D)$-triangularizable) if and
only if the invariant factors are Wedderburn polynomials
(\ref{diagonalization general case}) (resp. product of linear
polynomials (\ref{triangulation general case})). We also define
and study left and right eigenvalues of a matrix $A\in M_n(K)$ and
get analogues of classical results for commutative polynomials.
The last section is concerned with the notion of $G$-algebaric
sets. They give, in particular, another approach to the
Wedderburn's theorem on factorization of central polynomials. In
this last section we only consider the "classical" case i.e. we
assume that $S=id.$ and $D=0$.

\sk

\bs

\section {\bf Recapitulation}

Let us start with a brief review of basic definitions, notations
and contents of our previous paper "Wedderburn polynomials over
division rings, I".  We will refer this paper by "Wed1"(Cf.
\cite{LL5}).  Let us start with a triple $(K,S,D)$, where $K$ is a
division ring, $S$ is a ring endomorphism of $K$, and $D$ is a
$(S,Id.)$-derivation on $K$.  The latter means that $D$ is an
additive endomorphism of $K$ such that, for $a,b \in K ,\;
D(ab)=S(a)D(b)+D(a)b$. In the sequel of the paper a
$(S,\id)$-derivation will just be called a $S$-derivation.  We
will occasionally need the symmetric notion of a
$(\id,\sigma)$-derivation, $\delta$, where $\sigma$ is an
endomorphism of $K$ and $\delta$ is an additive map such that, for
$a,b\in K$, $\delta(ab)= a\delta(b)+\delta(a)\sigma(b)$.  In
particular, when $S$ is an automorphism of $K$ and $D$ is an
$S$-derivation, the map $-DS^{-1}$ is a $(Id.,S^{-1})$-derivation.

In the general $(K,S,D)$-setting, we can form the Ore ring of skew
polynomials $\,K[t;S,D]$. More details about this ring and its
properties can be found in the introduction of "Wed1" or in
\cite{Co3}.

In case $\,D=0 \,$ (resp. $S=I$), we write $\,K[t;S]\,$ (resp.
$K[t;D]$) for the skew polynomial ring $K[t;S,0]$ (resp. $K[t;Id.
,D]$). Of course, when $\,(S,D)= (Id.,0)\,$ (we refer to this as
the ``classical case''), $\,K[t;S,D]\,$ boils down to the usual
polynomial ring $\,K[t]\,$ with a {\it central\/} indeterminate
$\,t$.  {\it Throughout this paper, we'll write $\,R:=K[t;S,D]$}.
$R$ is a right euclidian domain (hence, in particular, a left
principal domain).   For $f(t) \in R$ and $a \in K$ there exist
$q(t)\in R$ and $b \in K$ such that
$$ f(t) = q(t)(t-a) + b \, , \; {\rm we \;then \; define} \; f(a):=b $$
For details Cf. \cite{LL1}, \cite{LL2} or Wed1. A subset $\De
\subseteq K$ is algebraic if there exists a polynomial $g \in R $
such that $g(x) =0$ for all $x\in \De$.  For $f \in R$ we put
$V(f):=\{ a \in K \vert f(a) = 0 \}$.  This set is obviously
algebraic and we say that a polynomial $f \in R$ is a Wedderburn
polynomial if $f$ is monic and is of minimal degree amongst
polynomials annihilating $V(f)$.  An element $a \in K$ is
$P$-dependent over an algebraic subset $\De$ if any polynomial
annihilating $\De$ also annihilates $a$.  A subset $B$ of an
algebraic set $\De$ is called a $P$-basis for $\De$ if no element
$b \in B$ is $P$-dependent over $B\setminus \{b\}$ and all
elements of $\De$ are $P$-dependent over $B$.  The cardinal of a
$P$-basis is called the rank of the algebraic set and is denoted
$\rk \De$.

An element $b \in K $ is $(S,D)$-conjugate to an element $a\in K$
if there exists $c\in K \setminus \{0\}$ such that $b= S(c)ac^{-1}
+ D(c)c^{-1}$, in this case we write $ b:=a^c $ and the set $\{a^x
\vert x\in K\setminus \{0\} \}$ will be denoted $\D^{S,D}(a)$ (or
just $\D (a)$ when no confusion is possible) and called the
$(S,D)$-conjugacy class of $a$. For $a\in K$ we define the
$(S,D)$-centralizer of $a$, denoted by $C^{S,D}(a)$, to be the set
$C^{S,D}(a):=\{x\in K \setminus \{0\}\; \vert \; a^x=a \}\cup
\{0\} $. This is in fact a division subring of $K$. Of course
these notions have analogues for the case of a
$(\id,\sigma)$-derivation $\delta$. For instance an element $b\in
K$ is $(\delta,\sigma)$-conjugate to an element $a\in K$ if there
exists $c\in K \setminus \{0\}$ such that $b=ca\sigma(c^{-1}) +
c\delta(c^{-1})$. The set of elements $(\delta,\sigma)$-conjugate
to an element $a$ will be denoted $\D^{\delta,\sigma}(a)$.  It is
an easy exercise to remark that, when $\sigma$ is an automorphism
of $K$, we have $\D^{S,D}(a)=\D^{-DS^{-1},S^{-1}}(a)$ (Cf.
\ref{comparison left and right}).

For $h \in R$ and $x\in K \setminus V(h)$ we define
$\phi_h(x):=x^{h(x)}$. This map appears naturally while evaluating
a product $gh$ at an element $x\in K \setminus V(h)$:
\begin{equation}\label{eq1}
gh(x)= g(\phi_h(x))h(x).
\end{equation}

Let us recall that $\phi_h(\D(a))\subseteq \D(a)$ i.e.\,$\phi_h$
preserves the $(S,D)$-conjugacy classes. While computing $\phi_h$
within a single $(S,D)$-conjugacy class $\D (a)$, another map
naturally appears: $\lambda_{h,a}: K \longrightarrow K: x \mapsto
h(a^x)x$.  An easy exercise shows that, if $a^x \in K\setminus
V(h)$, we have $\phi_h(a^x)=a^{\lambda_{h,a}(x)}$. The map
$\lambda_{h,a}$ is in fact right $C:=C^{S,D}(a)$-linear and
$\ker\lambda_{h,a}=\{x\in K\setminus \{0\} \, \vert \, a^x \in
V(h)\}\cup \{0\}$.  Moreover if an algebraic set $\Gamma$ is
contained in a conjugacy class $\D(a)$, say $\Gamma=a^Y$ for some
$Y\subseteq K\setminus \{0\}$, then $V(f_\Gamma)=a^{YC}$, where
$YC$ is the right $C=C^{S,D}(a)$-vector space generated by $Y$ and
$\rk\Gamma=\deg f_\Gamma=\dim_CYC$ (Cf \cite{LL2}). We also have
$\rk (V(h)\cap \D (a))= \dim_C\ker\lambda_{h,a}$ (Cf Wed1). Let us
also remark that, for $f,g\in R$, we have
$\lambda_{fg,a}=\lambda_{f,a}\lambda_{g,a}$.

\section {\bf Rank theorems}

In this section we will present different relations involving the
rank of an algebraic set.   Our first objective is to relate the
rank of $V(gh)$ and the ranks of $V(g)$ and $V(h)$.  Let us first
recall the following result from Wed1 (Cf. \cite[Corollary
4.4]{LL5}).

\begin{lemma}
\label{decomposition} If $\,\D_i \;(1\leq i \leq r)\,$ are
algebraic sets located in different $\,(S,D)$-conjugacy classes
$\D^{S,D} (a_i)$ of $\,K$, then
\begin{enumerate}
\item The set $E_i:=\{x\in K\setminus \{0\} \, \vert \,
a_i^x \in \D_i\}\cup \{0\}$ is a right vector space over
$C_i:=C^{S,D}(a_i)$.
\item $$\,\rk\bigl(\,\bigcup_{i=1}^r
\D_i \bigr)= \sum_{i=1}^r\, \rk\D_i=\sum_{i=1}^r\,\dim_{C_i}E_i.
$$
\end{enumerate}

\end{lemma}

Of course, this lemma applies to the set $V(f)$ of right roots of
a polynomial $f\in R$.  For $f\in R=K[t;S,D]$ and $a\in K$, we
denote $V(f)=\{x\in K \, | \, f\in R(t-x)\}$, $V'(f)=\{x\in K\,
|\, f\in (t-x)R \}$, $E(f,a)=\{x\in K\setminus \{0\}\,|\, a^x\in
V(f)\}\cup \{0\}$.  $E(f,a)$ is a right $C^{S,D}(a)$-vector space.

\begin{corollary}
\label{Conjugacy classes intersecting roots}

 With the above notations one has:
\begin{enumerate}
\item $V(f)$ intersects at most $n=deg(f)$ $(S,D)$-conjugacy classes,
say $V(f)=\cup_{i=1}^r(V(f)\cap \Delta(a_i) )$, with $r\le n$.
\item $$
rk V(f)= \sum_{i=1}^rdim_{C_i}E(f,a_i)\le deg(f),\; {\rm where}\;
C_i=C^{S,D}(a_i).
$$
The equality holds if and only if $f$ is a Wedderburn polynomial.
\item $V'(f)\cup V(f)$ intersects at most $n=deg(f)$ $(S,D)$-conjugacy classes.
\end{enumerate}

\end{corollary}
\begin{proof}
(1).  Let us recall that any polynomial $f\in R=K[t;S,D]$ can be
factorized as a product of irreducible polynomials: $f=p_1\cdots
p_n$.   Moreover if $f=q_1\cdots q_l$ is another such
factorization then $l=n$ and there exists a permutation $\pi\in
S_n$ such $R/Rp_i\cong R/Rq_{\pi(i)}$ (this means that $R$ is a
UFD, Cf. \cite{Co2}).On the other hand, it is easy to check that
$R/R(t-a)\cong R/R(t-b)$ if and only if $\Delta(a)=\Delta(b)$ (see
Thm. \ref{similarity of polynomilas and (S,D)conjugation} for a
further generalization).  It is then clear that the number of
conjugacy classes containing right roots of $f$ is bounded by
$deg(f)$.

Alternatively one can apply the above lemma \ref{decomposition} to
the algebraic set $V(f)$ to prove this result.  This is left to
the reader.

\noindent (2). Decomposing $V(f)$ into the $(S,D)$-conjugacy
classes it intersects, we can write $V(f)=\cup_{i=1}^r\Delta_i$
where $\Delta_i=V(f)\cap\Delta(a_i)$ and $E(f,a_i)=\{x\in
K\setminus \{0\} \,|\, f(a_i^x)=0\}\cup \{0\}$.  The above lemma
\ref{decomposition} then yields the desired formulas and the
additional statement comes from the fact that $f$ is a Wedderburn
polynomial if and only if $rk(V(f)) =deg (f)$.

\noindent (3).  As in (1) above, this is again a direct
consequence of the fact that $R=K[t;S,D]$ is a UFD.
\end{proof}

Notice the following important special case: $E(f,0)$ is easily
seen to be the solution space of the differential equation
$f(D)=0$ and $C^{S,D}(0)=K_D$ is the constant subdivision ring of
$K$.  Amitsur's well-known theorem states that the dimension over
$K_D$ of the solution space of the equation $f(D)=0$
 is bounded by the degree of the polynomial $ f$.  This is now clear:
this dimension is one of the dimension appearing in the expression
of $\rk V(f)$.

\begin{lemma}
\label{linear algebra} Let $V$ be a right vector space over a
division ring $C$ and $\phi,\psi \in End_CV$.  If
$v_1,v_2,\dots,v_r$ is a basis for $\ker \psi$ and
$u_1,u_2,\dots,u_s \in V \setminus \ker\psi$ the following are
equivalent:
\begin{enumerate}
\item[i)] The set $\{v_1,\dots,v_r,u_1,\dots,u_s\}$ is a basis
for $\ker \phi \psi$.
\item[ii)] The set $\{\psi(u_1),\psi(u_2),\dots,\psi(u_s)\}$ is a
basis for $\im\psi \cap \ker\phi$.
\end{enumerate}
In particular, we have

$$ \dim_C\ker(\phi \psi) = \dim_C\ker\psi + \dim_C(\im\psi \cap
\ker\phi).$$
\end{lemma}
\begin{proof}
The easy proof is left to the reader as an exercise in linear
algebra.
\end{proof}

\ms

\begin{theorem}
\label{first rk thm} Let $g,h$ be polynomials in $R$, then
$$\rk V(gh)= \rk V(h) + \rk (\im\phi_h \cap V(g)). $$
In particular, we always have $$\rk V(gh) \le \rk V(h) + \rk
V(g).$$

\end{theorem}
\begin{proof}
Let us put $f=gh$ and remark that, thanks to Lemma
\ref{decomposition}, it is enough to prove that, for any $a \in
K$, we have $\rk(V(gh) \cap \Delta (a)) = \rk(V(h) \cap \Delta
(a)) + \rk(\im\phi_h \cap V(g) \cap \Delta(a)) $.  Using the
definitions and results recalled at the end of section $2$, we
get, for $a$ in $K$, $\lambda_{f,a} = \lambda_{g,a}\lambda_{h,a}$.
In particular, $\ker\lambda_{h,a} \subseteq \ker\lambda_{f,a}$.
Moreover, if $C$ stands for $C^{S,D}(a)$, we have $\rk(V(f) \cap
\Delta (a)) = \dim_C\ker\lambda_{f,a}$ ; $\rk(V(h) \cap \Delta(a))
= \dim_C\ker\lambda_{h,a}$ ; $\im\phi_h \cap \Delta(a) =
a^{\im\lambda_{h,a}\setminus\{0\}}$ and $\rk(V(g) \cap \im \phi_h
\cap \Delta(a)) = \dim_C(\im\lambda_{h,a} \cap
\ker\lambda_{g,a})$. So we finally must prove that

$$\dim_C\ker\lambda_{f,a} = \dim_C\ker\lambda_{h,a} +
\dim_C(\im\lambda_{h,a} \cap \ker\lambda_{g,a}).$$

But this is exactly what is given by Lemma \ref{linear algebra}.
\end{proof}

As an application of the above result let us give another proof of
the main part of the "factor theorem" \cite{LL5}, Theorem 5.1.
Recall that $f\in \mathcal W $ if and only if $f$ is monic and
$\rk V(f) = \deg f$
\begin{corollary}
\label{factors of W polynomials} If $ f = gh \in \mathcal W$ then
$g,h \in \mathcal W$
\end{corollary}
\begin{proof}
  The above theorem
implies that $\rk V(g) + \rk V(h) \ge \rk V(gh)=\deg f = \deg g
+\deg h$.  This implies $\rk V(g) = \deg g$ and $\rk V(h) = \deg h
$.
\end{proof}

Recall from Wed1, that if $\Delta \subseteq K$ is an algebraic
set, we denote by $f_\Delta$ the monic polynomial of minimal
degree annihilating $\Delta$, and we put $\overline \Delta = \{x
\in K \vert \, f_\Delta(x)=0 \}$.
\begin{theorem}
\label{second rk thm} Let $h \in R$ and $\Delta \subseteq K$ be an
algebraic set disjoint from $V(h)$.  Then:
\begin{enumerate}
\item $\phi_h(\D)$ is an algebraic set.
\item $$\rk \phi_h(\Delta) = \rk \Delta - \rk(\overline \Delta \cap V(h)).$$
\item $\rk\phi_h(\Delta) = \rk\Delta$ iff $\overline \Delta \cap V(h)
= \emptyset$.
\end{enumerate}
\end{theorem}
\begin{proof}
1. Let $g,g'\in R$ be such that $[f_\D,h]_l=gh=g'f_\D$.  Then for
$x\in \D$, we have $0=(g'f_\D)(x)=(gh)(x)=g(\phi_h(x))h(x)$.
Hence, since $h(x)\ne 0$, $\phi_h(x)=0$.

\noindent 2. Decomposing the algebraic sets $\Delta,
\phi_h(\Delta)$ and $\overline \Delta \cap V(h)$ in conjugacy
classes and using the above lemma \ref{decomposition} we see that
it is enough to show that, for any $a\in K, \;\; \rk
(\phi_h(\Delta) \cap \Delta(a)) = \rk(\Delta \cap \Delta(a)) -
\rk(\overline \Delta \cap V(h)\cap \Delta(a))$.  Put $Y:=\{y \in K
\setminus \{0\} \vert a^y \in \Delta \cap \Delta(a) \}$ and denote
by YC the right $C^{S,D}(a)$-space generated by $Y$. We have $\rk
(\Delta \cap \Delta(a))= \rk(\{a^y \vert y \in Y\}) = \dim_CYC \;
; \; \rk(\overline \Delta \cap V(h) \cap \Delta(a)) = \rk( \{a^y
\vert y \in YC \; {\rm and} \; h(a^y) = 0\} )= \dim_C(YC\cap
\ker\lambda_{h,a})$ and $\rk(\phi_h(\Delta) \cap \Delta(a)) =
\dim_C\lambda_{h,a}(YC)$.  Consider  the map $\lambda_{h,a}$
restricted to $YC$ ; the required equality is an immediate
consequence of the relation between the dimension of the kernel
and the dimension of the image of this map.

\noindent 3. This is a particular case of 2. above.
\end{proof}

\begin{example}
{\rm Let $K$ be a division ring (we assume that $S=id.,D=0$) and
$a,x\in K,\;x\notin\{0,-1\}$, be such that $\{a,a^x,a^{1+x}\}$ are
distinct elements. Consider the polynomial $h(t)=t-a^{1+x}\in
K[t]$ and $\Delta = \{a,a^x\}$.  It is easy to check that
$V(h)\cap \Delta = \emptyset,\; V(h)\cap \overline \Delta =
\{a^{1+x}\}$.  Notice also that
$h(a^x)x=a^xx-(1+x)a+a^{1+x}=-a+a^{1+x}=-h(a)$ and thus
$\phi_h(a^x)=a^{h(a^x)x}=a^{h(a)}=\phi_h(a)$.  This gives
$\phi_h(\Delta)=\{a^{a-a^{1+x}}\}$. Of course, the above formula
can be checked on this particular example.  This also shows that
it is necessary to take $V(h) \cap \overline \Delta$ and not
merely $V(h) \cap \Delta$ in the formula.}
\end{example}

As a corollary let us mention the following interesting fact:

\begin{corollary}
For $h\in R$, let $\{a_1,\dots,a_n\}$ be a $P$-basis for $V(h)$
and $\{b_1,\dots,b_s\}\subset K \setminus V(h)$.  Then
$\{a_1,\dots,a_n,b_1,\dots,b_s\}$ is $P$-independent if and only
if $\{\phi_h(b_1),\dots,\phi_h(b_s)\}$ is $P$-independent.
\end{corollary}
\begin{proof}
The proof follows easily from the above theorem if we put
$\De=\{b_1,\dots,b_s\}$ and remark that $\overline{\De} \cap V(h)
= \emptyset$ iff $\{a_1,\dots,a_n,b_1,\dots,b_s\}$ is
$P$-independent.
\end{proof}

\section{\bf Companion matrices}

In this section we will show that the companion matrices together
with pseudo linear transformations give a natural interpretation
of some notions related to $R=K[t;S,D]$-modules.

\begin{definition}
 Two polynomials $g,h\in R=K[t;S,D]$ are similar if $R/Rg \cong
R/Rh$.   This will be denoted by $f\sim g$.  $\Delta(f)$ will
stand for the set of polynomials similar to $f$.
\end{definition}

\begin{remark}
{\rm
 The notion of similarity can be introduced over a general ring.
It is obviously an equivalence relation and in an integral domain
we always have $R/Rg\cong R/Rf$ if and only if $R/gR\cong R/fR$
(Cf. \cite{LO}, \cite{Co2}).
 }
\end{remark}

\begin{example}
{\rm  Let $a,b\in K$, then $t-a\sim t-b$ if and only if $a$ and
$b$ are $(S,D)$-conjugate.  Theorem \ref{similarity of polynomilas
and (S,D)conjugation} will generalize this example and give a
description of similarity of polynomials in terms of
$(S,D)$-conjugation.}
\end{example}

\begin{lemma}
\label{left common multiple} Let $f,g,h\in R$ be monic
polynomials. Then:
\begin{enumerate}
\item There exist uniquely determined monic polynomials $g',h' \in R$
such that $Rg \cap Rh =Rg'h=Rh'g$. We will denote $g'$ and $h'$ by
$g^h$ and $h^g$ respectively.
\item
$$
\frac{R}{Rh^g} \cong \frac{Rg+Rh}{Rh}
$$
In particular if $Rg +Rh=R$ we have $h^g\sim h$ and hence $\deg
h^g=\deg h$.
\item
$$
 Rfg \cap Rh =
\begin{cases}
Rfg & \text{if}\; g \in Rh, \\
(Rf \cap Rh^g)g & \text{if}\; g \notin Rh.
\end{cases}
$$
\end{enumerate}
\end{lemma}
\begin{proof}
1) This is clear.

2) This is given by a classical isomorphism theorem.  Notice also
that the map $R/Rh^g \longrightarrow R/Rh: x \mapsto xg$ is easily
seen to be well defined and injective. Moreover it is onto when
$Rg+Rh=R$.

3) This is easy to check and is left to the reader.

\end{proof}

\begin{remark}
{\rm Let us first notice that if $g=t-a$ and $h=t-b, \, a\ne b,$
we have $g^h=t-a^{a-b}$, where, as usual,
$a^c=S(c)ac^{-1}+D(c)c^{-1}$ for $c\in K\setminus\{0\}$.  More
generally, when $h=t-a$ we have $Rg \cap R(t-a)=Rg$ if $g(a)=0$
and $Rg \cap R(t-a)=R(t-a^{g(a)})g$ if $g(a)\ne 0$.  Remark also
that, when $h=t-a$, the formula in \ref{left common multiple}($3$)
above gives back the way of evaluating the product $fg$ at the
element $a\in K$. }
\end{remark}

We collect without proofs some easy facts related to similarity.
\begin{lemma}
\label{properties of similarity}

For $f,g,h$ monic polynomials in $R$ we have:
\begin{enumerate}
\item $\deg f^g\le \deg f$.
\item If $g-h\in Rf$, then $f^g\sim f^h$.
\item $\Delta (f)=\{f^q \, \vert \, q\in R, \; Rq + Rf =R \; {\rm
and}\; \deg q < \deg f\}$.
\item $gh\in Rf$ if and only if either $h\in Rf$ or $g\in Rf^r$
where $r$ is the remainder of $h$ right divided by $f$.
\item $(f^g)^h=f^{hg}$.
\end{enumerate}
\end{lemma}
\begin{proof}
We leave the proofs of these statements to the reader (Cf
\cite{LO} for similar facts in the more general frame of
$2$-firs).
\end{proof}

For a monic polynomial $f(t) = \sum_{i=0}^na_it^i \in R=K[t;S,D]$,
the companion matrix of $f$ denoted by $C_f$ is the $n \times n$
matrix defined by

$$ C_f = \begin{pmatrix}0&1&0&\cdots&0 \\0&0&1&\cdots&0 \\
\vdots&\vdots&\vdots&\ddots&\vdots \\
0&0&0&\cdots&1\\
-a_0&-a_1&-a_2&\cdots&-a_{n-1}
\end{pmatrix}.$$
We need also some results on pseudo-linear transformations
(abbreviated $PLT$ or $(S,D)$-PLT in the sequel).  For details on
this topic we refer the reader to \cite{L}, for instance. Let us
recall that for a left $K$-vector space $V$, a map $T: V
\longrightarrow V$ is an ($S,D$)-PLT if $T$ is additive and
$T(\alpha v)=S(\alpha)T(v)+D(\alpha)v$ for $\alpha \in K$ and
$v\in V$.  Let $A$ be a matrix in $M_n(K)$ and let $K^n$ stand for
the set of row vectors with coefficients in $K$. The maps $S$ and
$D$ can be extended to $K^n$ and to $M_n(K)$ in the obvious way.
Define the map $T_A:K^n \longrightarrow K^n:v \mapsto S(v)A +
D(v)$.  $T_A$ is an ($S,D$)-PLT which defines a left
$R=K[t;S,D]$-module structure on $K^n$ via
$(\sum_{i=0}^n\alpha_it^i).v=\sum_{i=0}^n\alpha_i(T_A)^i(v)$ for
$v\in K^n$ and $\sum\alpha_it^i\in K[t;S,D]$.  Conversely any
structure of left $R$-module defined on $K^n$ is of this form. Let
us denote $e_i:=(0,\dots,1,0\dots O)$ the element of $K^n$ with a
one in position $i$ and zero elsewhere.  For a monic polynomial
$f\in R$ of degree $n$, the $K$-linear map $R/Rf\longrightarrow
K^n:t^i \mapsto e_{i+1}, {\rm for} i=0,1,\dots,n-1$ induces an
$R$-module structure on $K^n$ that corresponds to $T_{C_f}$ where
$C_f$ is the companion matrix defined above.  The matrix
representing a $PLT$ depends on the $K$-basis of $K^n$ which is
chosen.  If two matrices $A$ and $B$ represent the same $PLT$ in
different bases, there exists an invertible matrix $P \in GL_n(K)$
such that
$$
B:= S(P)AP^{-1} + D(P)P^{-1}.
$$

This leads to the following definitions.

\begin{definitions}
\label{diagonalization, triangulation}
\begin{enumerate}
\item Two matrices $A,B \in M_n(K)$ are $(S,D)$-similar if there
exists an invertible matrix $P \in GL_n(K)$ such that $B=
S(P)AP^{-1} + D(P)P^{-1}$.
\item A matrix $A$ is $(S,D)$-diagonalizable (resp.
triangularizable) if it is $(S,D)$-similar to a diagonal (resp.
triangular) matrix.
\end{enumerate}
\end{definitions}

\begin{lemma}
\label{companion matrix of a product} Let $f \in R=K[t;S,D]$ be a
monic polynomial. Then:
\begin{enumerate}
\item All submodules of $R/Rf$ are of the form $Rg/Rf$, where $g$
is a monic right factor of $f$.
\item If there exist $a_1,\dots
,a_n \in K$ such that $f(t)=(t-a_n)(t-a_{n-1})\cdots (t-a_1)$,
then the companion matrix $C_f$ is $(S,D)$-similar to the
following one:
$$
\begin{pmatrix}
a_1     &      1 &        0 &       0&   \cdots &0        &0\\
  0     &    a_2 &        1 &       0&   \cdots &0        &0\\
  0     &      0 &       a_3&       1&          &         &  \\
\vdots  &        &   \ddots & \ddots &  \ddots  &         &\vdots \\
  0     &        &          &        &          & 1       &0\\
  0     &        &          &        &          & a_{n-1} & 1 \\
  0     &        &          &  \cdots&          &    0    & a_n \\
\end{pmatrix}
$$
\item If $f=gh$ where $g,h \in R$ are monic then the companion
matrix $C_f$ is $(S,D)$-similar to the following matrix
$$
\begin{pmatrix}
C_{h} & \begin{matrix}
        0 & \cdots & 0 \\
        \vdots & \cdots & \vdots \\
        1 & \cdots & 0 \\
      \end{matrix} \\
&    &    \\
 0  & C_{g} \\
\end{pmatrix}
$$
Where the rectangular matrices are of the required sizes.
\end{enumerate}
\end{lemma}
\begin{proof}
1) This is clear since $R$ is a left principal domain.

2) Notice first that the set $\{1 + Rf, t - a_1 + Rf,
(t-a_2)(t-a_1) + Rf, \dots , (t-a_{n-1})(t-a_{n-2})\cdots(t-a_1) +
Rf \}\subseteq R/Rf$ is a $K$-basis of $R/Rf$.  In this $K$-basis
the matrix associated to
 left multiplication by $t$ on $R/Rf$ is exactly the one
 displayed in the statement $2)$.  This shows that $C_f$ is
 $(S,D)$-similar to this matrix.

3) Put $l= \deg g$ and $n = \deg h$.  It is enough to consider the
following $K$ basis of $R/Rf$:

$ 1 + Rf, t +Rf , \dots, t^{n-1} +Rf , h +Rf , th +Rf, \dots ,t^{l
-1}h +Rf. $

It is easy to check that in this basis the matrix representing
left multiplication by $t$ is exactly the one mentioned in the
statement of the lemma.  This shows that this matrix is
$(S,D)$-similar to $C_f$.
\end{proof}

Let us remark that the second statement in the above lemma
\ref{companion matrix of a product} could also be obtained by
using the third one repeatedly.

The following easy lemma will be very useful allowing us to
translate $R=K[t;S,D]$-module theoretic notions into matrix
related ones. It will be used again in the next section.

\begin{lemma}
\label{matrix of a left R-morphism} Let $_RV$ and $_RW$ be left
$R$-modules which are finitedimensional as left $K$-vector spaces
with bases $\mathcal B$ and $\mathcal C$ respectively. Let
$\varphi:V \longrightarrow W $ be a left $K$-linear map and denote
$$
P:=M^{\mathcal B}_\mathcal C(\varphi) \quad A:=M^\mathcal
B_\mathcal B(t.) \quad {\rm and }\quad B:=M^\mathcal C_\mathcal
C(t.)\, .
$$
Then $\varphi$ is a morphism of left $R$-modules if and only if
$AP =S(P)B + D(P)$.
\end{lemma}

\begin{proof}
For a vector $v\in V$ we denote $v_{\mathcal B}$ the row in $K^n$
consisting of the coordinates of $v$ in the basis $\mathcal B$. We
use similar notations in $W$.  The definition of $M^\mathcal
B_\mathcal B(t.)$ gives that $(t.v)_\mathcal B =S(v_\mathcal B)A +
D(v_\mathcal B)$ and so $\varphi (t.v)_\mathcal C=S(v_\mathcal
B)AP +D(v_\mathcal B)P$. On the other hand,
$(t.\varphi(v))_\mathcal C =S(\varphi(v)_\mathcal
C)B+D(\varphi(v)_\mathcal C)=S(v_\mathcal B P)B + D(v_\mathcal
BP)=S(v_\mathcal B)(S(P)B+D(P))+D(v_\mathcal B )P$. Since
$\varphi$ is a morphism of left $R$-modules if and only if
$\varphi\circ t.=t.\circ \varphi$, we obtain the required
equality.
\end{proof}

As a first consequence we get the following:

\begin{theorem}
\label{similarity of polynomilas and (S,D)conjugation} Two monic
polynomials $f,g\in R$ are similar if and only if their companion
matrices $C_f$ and $C_g$ are $(S,D)$-conjugate.
\end{theorem}
\begin{proof}
Let $\mathcal B:=\{1+Rf,t+Rf,\dots,t^{n-1}+Rf\}$, where $n=\deg
f$, be a basis for the left $K$-vector space $R/Rf$.  Then $C_f$
represents the $(S,D)$-pseudo linear transformation $t.$ acting on
$R/Rf$ i.e. $C_f=M^{\mathcal B}_{\mathcal B}(t.)$. Similarly $C_g$
represents $t.$ in the appropriate basis $\mathcal C$ of $R/Rg$.
Since $f\sim g$ if and only if there exists an isomorphism $R/Rf
\stackrel{\varphi} \cong R/Rg $ of left $R$-modules.  Hence the
matrix $P:=M^{\mathcal B}_{\mathcal C}(\varphi)$ is invertible and
the above lemma \ref{matrix of a left R-morphism} shows $f\sim g$
that $C_f$ and $C_g$ are $(S,D)$-conjugate.

\end{proof}

\begin{prop}
\label{1 in Rg + hR and the metro equation} Let $g,h \in
R=K[t;S,D]$ be two monic polynomials of degree $l$ and $n$
respectively.  Put
$$
A:= \begin{pmatrix}
C_h & U \\
0 & C_g \\
\end{pmatrix}
\quad  {\rm and } \quad B:= \begin{pmatrix}
C_h & 0 \\
0 & C_g \\
\end{pmatrix}
$$
where $C_g,C_h$ denote the companion matrices of $g$ and $h$
respectively and $U$ is the unit matrix $e_{n1}\in M_{n\times
l}(K)$. Then the following are equivalent:
\begin{enumerate}
\item $0 \longrightarrow R/Rg {\stackrel{.h}\longrightarrow}R/Rgh
\longrightarrow R/Rh \longrightarrow 0 $ splits.
\item $1 \in Rg + hR$.
\item There exists a matrix $X\in M_{n\times l}(K)$ such that
$$
\begin{pmatrix}
I & S(X) \\
0 & I \\
\end{pmatrix}
A + \begin{pmatrix}
0 & D(X) \\
0 & 0 \\
\end{pmatrix}=
B
\begin{pmatrix}
I & X \\
0 & I \\
\end{pmatrix}
$$

\item There exists a matrix $X\in M_{n \times l}(K)$ such that
$C_hX-S(X)C_g-D(X)=U$ where $U$ is the matrix unit $e_{nl}$.
\end{enumerate}
\end{prop}
\begin{proof}
$(1)\Rightarrow (2)$ By hypothesis there exists a map $\varphi
:R/Rgh \longrightarrow R/Rg$ such that $\varphi\circ
.h=id._{R/Rg}$. Let $y\in R$ be such that $\varphi(1+Rgh)=y+Rg$.
We then have $(\varphi\circ .h)(1+Rg)=1+Rg$, i.e. $hy-1\in Rg$.
This gives that there exists $x\in R$ such that $hy+xg=1$.

%$(2)\Rightarrow (1)$ If $1=xg+hy$, we define $\varphi
%:R/Rgh\longrightarrow R/Rh\bigoplus R/Rg: u + Rgh \mapsto
%(u+Rh,uy+Rg)$.  It is easy to check that this map is a well
%defined isomorphism of left $R$-modules and that the map
%$p\circ\varphi$, where $p$ is the projection $R/Rh \bigoplus R/Rg
%\longrightarrow R/Rg $, splits the map $.h$.

$(2)\Rightarrow (3)$ By hypothesis there exist $x,y\in R$ such
that $1=xg+hy$.  using the right euclidian division, we may assume
that $deg(y)<deg(g)$.  Define $\varphi: R/Rgh \longrightarrow R/Rh
\oplus R/Rg: u+Rgh \mapsto (u+Rh,uy+Rg)$. It is easy to check that
this map is a well defined morphism of left $R$-modules.  Let
$\mathcal B = \{1+Rgh,t+Rgh,\dots,t^{n-1}+Rgh,h+Rgh,th+Rgh,\dots,
t^{l-1}h +Rgh \}$ and $\mathcal
C:=\{(1+Rh,0),(t+Rh,0),\dots,(t^{n-1}+Rh,0),(0,1+Rg),\dots,(0,t^{l-1}+Rg)\}$
be bases for the left $K$-vector spaces $R/Rgh$ and $R/Rh\bigoplus
R/Rg$, respectively. Since $hy+Rg=1+Rg$, it is easy to check that
the matrix of $\varphi$ in these bases is of the form
$$
P:=M^\mathcal B_\mathcal C(\varphi)=
\begin{pmatrix}
I & Y \\
0 & I \\
\end{pmatrix}
$$
(where $Y$ is the $n\times l$ matrix whose rows are given by
writing $t^iy + Rg$, $i=1,\dots ,n-1$, in the basis $t^j + Rg$,
$j\in \{0,\dots,l-1\}$).  Remark that we also have $A=M^{\mathcal
B}_\mathcal B(t.)$ and $B=M^\mathcal  C_\mathcal C(t.)$.  Since
$\varphi$ is a morphism of left $R$-modules, Lemma \ref{matrix of
a left R-morphism} implies that $AP=S(P)B+D(P)$ i.e. $S(P^{-1})A +
D(P^{-1})=BP^{-1}$.  We then get the desired conclusion with
$X:=-Y$.

 $(3) \Rightarrow (1)$ Let $\mathcal B$ and $\mathcal C$ be the bases for $R/Rgh$ and
$R/Rh\bigoplus R/Rg$ defined in the proof of $(2)\Rightarrow (3)$.
 Let $\varphi: R/Rgh \longrightarrow R/Rh \bigoplus R/Rg$ be the
left $K$-isomorphism map such that
$$P:=M^\mathcal B_\mathcal C(\varphi)=
\begin{pmatrix}
I & -X \\
0 & I \\
\end{pmatrix}
$$
We have $A=M^\mathcal B_\mathcal B(t.)$ and $B=M^\mathcal
C_\mathcal C(t.)$.  Statement $(3)$ implies that $S(P^{-1})A +
D(P^{-1})=BP^{-1}$ i.e. $AP=S(P)B + D(P)$.  The previous lemma
shows that $\varphi$ is in fact an homomorphism of left
$R$-modules.  Let $p$ denotes the projection $R/Rh\bigoplus R/Rg
\longrightarrow R/Rg$.  We claim that $p\circ\varphi:R/Rgh
\longrightarrow R/Rg$ is a splitting of $.h$. Indeed
$(p\circ\varphi\circ
.h)(1+Rg)=p(\varphi(h+Rgh))=p((0,1+Rg))=1+Rg$.

$(3)\Leftrightarrow (4)$ This is left to the reader.

\end{proof}

\section{diagonalization and triangulation}

In this section we will briefly consider a generalization of
Wedderburn polynomials called fully reducible polynomials. The
family of fully reducible polynomial is larger than the Wedderburn
one, but they share many properties and, for what we have in mind,
they are not more difficult to handle. They will show better the
connection between factorization in $R$ and companion matrices.
They were introduced by Ore himself and further studied by PM Cohn
in the setting of $2$-firs (\cite{Co2}) and more recently by the
second and third authors of this paper (again in the setting of
$2$-firs, Cf \cite{LO}).  The companion matrices of these families
of polynomials will lead us naturally to a characterization of
diagonalizability of a matrix over a division ring.

\begin{definition}
A monic polynomial $f\in R=K[t;S,D]$ is fully reducible if there
exist irreducible polynomials $p_1,\dots,p_n$ such that
$Rf=\cap_{i=1}^nRp_i$.
\end{definition}

 Wedderburn polynomials and monic irreducible polynomials are
fully reducible.  Notice also that a polynomial
$g(t)=(t-a_1)\cdots (t-a_n)$ is fully reducible if and only if it
is Wedderburn.

The notion of fully reducible polynomials is symmetric i.e. if $f
\in R=K[t;S,D]$ and $p_1,p_2,\dots,p_n$ are irreducible
polynomials such that $Rf=\cap_{i=1}^nRp_i$ then there exist
irreducible polynomials $q_1,\dots,q_n$ such that
$fR=\cap_{i=1}^nq_iR$.  Moreover there exists a permutation
$\pi\in S_n$ such that $p_i\sim q_{\pi(i)}$ i.e. $R/Rp_i\cong
R/Rq_{\pi(i)}$ (Cf. \cite{LL4} or \cite{LO}).

\begin{theorem}
\label{equivalent conditions for being fully reducible} Let $f\in
R$ be a monic polynomial of degree $l$.  Then the following are
equivalent:
\begin{enumerate}
\item $f$ is fully reducible.
\item There exist monic irreducible polynomials $p_1,\dots
,p_n$ such that $Rf=\cap_{i=1}^{n}Rp_i$ is an irredundant
intersection.
\item There exist monic irreducible polynomials $p_1,\dots,p_n\in R$
such that the map $\varphi: R/Rf \longrightarrow
\oplus_{i=1}^nR/Rp_i: q+Rf \mapsto (q+Rp_1,\dots,q+Rp_n)$ is an
isomorphism of $R$-modules.
\item There exist monic irreducible polynomials $p_1,\dots,p_n \in R$
and an invertible matrix $V\in M_l(K)$ such that
$$
C_fV=S(V)\diag(C_{p_{_1}},\dots,C_{p_n}) + D(V).
$$
\item $R/Rf$ is semisimple.
\end{enumerate}
\end{theorem}
\begin{proof}

$(1)\Leftrightarrow (2)$ is clear by definition.

$2) \Rightarrow 3)$.  The map $\varphi$ is is easily seen to be
well defined and injective. Since, for every $j\in\{1,\dots,n\}$,
$Rp_j+(\cap_{i\ne j}Rp_i)=R$, Lemma \ref{left common multiple}
shows that $\deg f=\sum_{i=1}^n\deg p_i$.  This implies that the
$\dim_K(R/Rf)=\dim_K(\oplus_i R/Rp_i)$ and we conclude that
$\varphi$ is onto.

$(3)\Longrightarrow (2)$.  Composing $\varphi$ with the natural
homomorphism $R{\stackrel{p}\longrightarrow}  R/Rf$ we obtain an
onto $R$-morphism: $\psi = \phi \circ p$ such that $Ker \psi= Rf$
and we conclude that $Rf=\cap_{i=1}^nRp_i$.  The fact the this
intersection is irredundant is clear from the equalities:
$l=deg(f)= dim_K(R/Rf)=\sum_i dim_K(R/Rp_i)=\sum_i deg(p_i)$.

$(3)\Rightarrow (4)$.  Let $\mathcal B=\{t^i+Rf \,\vert\,
i=0,\dots,l-1 \}$ be a basis for the left $K$ space $R/Rf$ and
 $\mathcal C=\{(0,\dots,0,t^j+Rp_i,0,\dots,0) \,\vert \, i=1,\dots,n \ {\rm and}\
 j=0,\dots,n_i-1\}$, where $n_i=\deg p_i$, be a $K$-basis for $\oplus_iR/Rp_i$.
 We have $M_{\mathcal B}^{\mathcal B}(t.)=C_f$ and
 $M_{\mathcal C}^{\mathcal C}(t.)=\diag (C_{p_1},\dots,C_{p_n})$.
 Put $V:=M_{\mathcal C}^{\mathcal B}(\varphi)$.  Then $V$ is invertible and since
 $\varphi$ is a morphism of left $R$-modules, lemma \ref{matrix of a left
 R-morphism} yields the required equality.

 $(4)\Rightarrow(3)$.  It is enough to define the map $\varphi$ via
 $M_{\mathcal C}^{\mathcal B}(\varphi)$ where $\mathcal B$ and $\mathcal
 C$ are the bases defined above.

 $(3)\Leftrightarrow (5)$.  This is clear and left to the reader.

\end{proof}

In (\cite{LL5}) (resp. \cite{LO}) several criterion were given for
a product of Wedderburn polynomials (resp. fully reducible
polynomials) to be again a Wedderburn polynomial (resp. fully
reducible). We will give two more criterions in the following
theorem.   We treat the cases of Wedderburn polynomials and fully
reducible polynomials simultaneously. Let us first introduce a
technical notation: For a polynomial $g=p_r\cdots p_1$ with $\deg
p_i=n_i$ for $i=1,\dots,r$, we put:
$$
C_g(p_r,\dots, p_1)=
    \begin{pmatrix}
          C_{p_1} & U_1 & 0 & \cdots \\
          0 & \ddots & \ddots & 0 \\
          0 & \cdots  & C_{p_{r-1}} & U_{r-1} \\
          0   &   0 &    0     & C_{p_r} \\
    \end{pmatrix},
$$
where for $i=1,\dots,r-1$, the matrices $U_i\in M_{n_i\times
n_{i+1}}(K)$ have a one in the bottom left corner and zero
elsewhere.
 In particular, if $g(t)=(t-a_r)\cdots (t-a_1)$ the above matrix
 takes the simpler form
 $$C_g(a_r,\dots , a_1)=
   \begin{pmatrix}
          a_1 & 1 & 0 & \cdots \\
          0 & \ddots & \ddots & 0 \\
          0 & \cdots  & a_{r-1} & 1 \\
          0   &   0 &    0     & a_r \\
    \end{pmatrix}.
 $$

Notice that, according to Lemma \ref{companion matrix of a
product}, this matrix represents the pseudo linear transformation
$t.$ acting on $R/Rg$ and hence is $(S,D)$-similar to $C_g$.

\begin{theorem}
\label{product of fully reducible polynomials} Let $g,h$ be fully
reducible polynomials (resp. $W$-polynomials) in $R$ of degree $l$
and $n$ respectively. Then the following are equivalent:
\begin{enumerate}
\item $gh$ is a fully reducible (resp. $W$-) polynomial.
\item $1\in Rg + hR$.
\item There exists a matrix $X \in M_{n\times l}(K)$ such that
$$C_hX - S(X)C_g - D(X) = U\, ,$$ where $U=e_{n1}\in M_{n\times l}(K)$.
\item If $g=p_r\cdots p_1$ and $h=q_s\cdots q_1$
(resp. $g=(t-b_l)\cdots (t-b_1)$ and $h=(t-a_n)\cdots (t-a_1)$)
There exists $Y \in M_{n\times l}(K)$ such that
$$
C_h(q_s,\dots, q_1)Y - S(Y)C_g(p_r,\dots, p_1) - D(Y) = U .
$$
(resp. $$C_h(a_n,\dots,a_1)Y - S(Y)C_g(b_l,\dots, b_1)-D(Y)= U.
\mbox{\rm ) }$$
\end{enumerate}
\end{theorem}

\begin{proof}
$(1)\Leftrightarrow (2)$ This comes from the fact that $gh$ is
fully reducible if and only if  $R/Rgh$ is semisimple and hence
the short exact sequence from Equation \ref{1 in Rg + hR and the
metro equation} splits and this proposition shows that $1\in Rg +
hR$.

$(2)\Leftrightarrow (3)$ This is exactly equivalence
$(2)\Leftrightarrow (4) $ of \ref{1 in Rg + hR and the metro
equation}.

$(2)\Leftrightarrow (4)$ This is obtained similarly as above
making use of the 2 bases in $R/Rgh$ we have used in Proposition
\ref{1 in Rg + hR and the metro equation}.  We leave the details
for the reader.

%One can check that if $\mathcal
%B=\{1+Rgh,t+Rgh,\dots,t^{l_1-1},q_1,tq_1,\dots\}$ is a basis for
%$R/Rgh$ and $\mathcal C$ is the basis of $(\oplus_{i=1}^s
%R/Rq_i)\bigoplus (\oplus_{j=1}^r R/Rp_j)$ consisting of powers of
%$t$ in each components, we have $M_\mathcal B^{}$ and $M_{\mathcal
%C}^{\mathcal C}(t.)=\diag (C_h(q_s\dot q_1)) $. On the other hand
%there exist bases $\mathcal B$ of $R/Rh$ and $\mathcal C$ of
%$R/Rg$ such that $M_{\mathcal B}^{\mathcal B}(t.)$
\end{proof}

In our previous work Wed1 (\cite{LL5}) we have obtained a few
conditions for a product of two $W$-polynomials to be a
$W$-polynomial.  Let us point out that the advantage of the
characterization $(3)$ in the above theorem is that there is a
finite number of equations to check and that they are directly
available from the coefficients of $g$ and $h$ themselves.   The
characterization $(4)$ is also interesting if one knows in advance
a factorization of $f$ and $g$.

\begin{example}
{\rm Let $K=\mathbb Q(x)$ be the field of rational fractions in
$x$ over the rational and let $R$ be the Ore extension $R=\mathbb
Q(x)[t;id.,\frac{d}{dx}]$.  Using the above theorem it is easy to
show that, for any $q\in \mathbb Q(x)$ and for any $n\in \mathbb
N$, the polynomials $(t-q)^n\in R$ are $W$-polynomials.  To check
this, let us write $(t-q)^n=(t-q)^{n-1}(t-q)$ and
$U=(1,0,\dots,0)\in M_{1\times n-1}(\mathbb Q(x))$.  Part $(4)$ of
the theorem, with $g=(t-q)^{n-1}$ and $h=t-q$, shows that we have
to find $(y_1,\dots,y_{n-1}) \in \mathbb Q(x)^{n-1}$ such that:
$$
\left\lbrace
     \begin{array}{l}
       y_1q + D(y_1) - qy_1 +1 = 0 \\
       y_1 + y_2q +D(y_2) -qy_2= 0  \\
       y_2 + y_3q +D(y_3) -qy_3= 0  \\
       \vdots \\
       y_{n-2} + y_{n-1}q +D(y_{n-1})-qy_{n-1}=0\\
     \end{array}
\right.
$$
It is then easy to see that the sequence defined by
$y_i=(-1)^{i+1}\frac{x^i}{i!}$ ($i=1,\dots,n-1$) gives a solution
of the above system of equations.  We can thus conclude that for
any $n\in \mathbb N$ the polynomial $(t-q)^n\in R$ is a
$W$-polynomial. }
\end{example}

\begin{example}
{\rm Let $k$ be a commutative field of characteristic $0$, $D$ a
derivation ($S=Id.$) on $k$.  Kolchin (Cf. \cite{Ko}) showed that
there exists a field $U$ containing $k$ as a subfield and a
derivation $\overline{D}$ over $U$ extending $D$ such that the
equation
$$
p(x,\overline{D}(x),\dots,\overline{D}^{(n)}(x))=0, \quad
n\;\,{\rm arbitrary},
$$
has a solution $u\in U$ for all $p(X)\in
U[X_1,\dots,X_{n+1}]\setminus U$.  Since for any $v \in U$ the
polynomial $X_2 - v$ has a solution, $\overline{D}$ is onto. {\it
We claim that all monic polynomials of $R=U[t;\overline{D}]$ are
W-polynomials}. Let us first show the the irreducible polynomials
are of degree at most $1$.  Indeed, if $p(t)=\sum a_it^i\in R$ is
such that $deg (p(t) )>1$ it is easy to verify that the hypothesis
made on $U$ implies that there exists $v\in U$ such that
$p(v)=\sum a_iN_i(v)=0$ i.e. $t-v$ divides $p(t)$ on the right. It
follows that any monic polynomial $h(t)$ of degree $n$ can be
factorized in the form $h(t)=(t-a_n)\dots (t-a_1)$.  By induction
on the degree we need only show that if $h(t)$ is a W-polynomial
than $(t-b)h(t)$ is also a W-polynomial.  Once again using the
above theorem \ref{product of fully reducible polynomials}(4), we
must find $(y_1,\dots , y_n)\in U^n$ such that :
$$
\begin{pmatrix}
          a_1 & 1 & 0 & \cdots \\
          0 & \ddots & \ddots & 0 \\
          0 & \cdots  & a_{n-1} & 1 \\
          0   &   0 &    0     & a_n \\
    \end{pmatrix}
\begin{pmatrix}
y_1\\
y_2\\
\vdots \\
y_n\\
\end{pmatrix} - \begin{pmatrix}
                y_1\\
                y_2\\
                \vdots \\
                y_n\\
                \end{pmatrix} b - \begin{pmatrix}
                                   D(y_1)\\
                                   D(y_2)\\
                                   \vdots \\
                                   D(y_n)\\
                                   \end{pmatrix}= \begin{pmatrix}
                                                  0 \\
                                                  0 \\
                                                  \vdots \\
                                                  1 \\
                                                  \end{pmatrix}.
                                                  $$
In other words we have to solve (for $y_i$'s) the equations
$$
a_iy_i - y_ib -D(y_i)= u_i \quad {\rm for}\; 1\le i \le n\, ,
$$
where $u_i=-y_{i+1}$ for $1\le i \le n-1$ and $u_n=1$.  But
solving first for $y_n$ and then for $y_{n-1}$,... it is easy to
check that these equations all have solutions thanks to the
property of $U$. }
\end{example}

We now come to the diagonalization.  As is well known, a matrix $A
\in M_n(k)$ over a commutative field $k$ is diagonalizable if and
only if its minimal polynomial can be written as a product of
distinct linear polynomials in $k[t]$.  In other words the minimal
polynomial of $A$ must be a W-polynomial. In the next section we
will generalize this result and obtain a criterion for the
diagonalizability of a matrix with coefficients in a division
ring.   This will be developed in an "$(S,D)$" setting.

Let us recall some results and notations from \cite{LL1}.  For
$\{b_1,\dots,b_n \} \subset K $ we define the Vandermonde matrix:

$$ V_n(b_1,\dots,b_n) = \begin{pmatrix}  1&1&\cdots&1\\ b_1&b_2&
\cdots & b_n \\ N_2(b_1)&N_2(b_2)&\cdots&N_2(b_n)\\
\vdots&\vdots&\vdots&\vdots\\N_{n-1}(b_1)&N_{n-1}(b_2)&
\cdots&N_{n-1}(b_n)
\end{pmatrix}  $$
where, for $a\in K$ and $i\ge 0$, $N_i(a)$ denotes the evaluation
of $t^i$ at $a$.  Notice that one has $N_0(a)=1$ and, using the
product formula recalled in ($2.1$), one gets
$N_{i+1}(a)=(tt^i)(a)=\phi_{t^i}(a)t^i(a)=S(N_i(a))a + D(N_i(a))$.

Let us also remark that this matrix appeared already in an hidden
form in \ref{equivalent conditions for being fully reducible}.
Indeed if, in this theorem, $p_1=t-b_1,\dots,p_n=t-b_n$ the matrix
$V$ in Theorem \ref{equivalent conditions for being fully
reducible} (4)(Cf. also its proof) is exactly the above
Vandermonde matrix. This can be exploited to get the equivalence
between $(iii)$ and $(iv)$ in the following proposition.
\begin{lemma}
\label{vandermonde and independence} For $\Delta:=\{b_1,\dots,b_n
\} \subset K$ the following are equivalent
\begin{enumerate}
\item[$i)$] $\Delta:=\{b_1,\dots,b_n \}$ is P-independent.
\item[$ii)$] $\deg f_{\Delta} = n$.
\item[$iii)$] $Rf_\Delta =\cap_{i=1}^nR(t-b_i)$.
\item[$iv)$] The matrix $V_n(b_1,\dots,b_n)$ is invertible.
\end{enumerate}
\end{lemma}
\begin{proof}
 i) $\Leftrightarrow$ ii)  and  ii) $\Leftrightarrow$ iii)  are easy
to establish and were proved in \cite{LL4},\cite{LL5}.

 (iii) $\Leftrightarrow$ iv)  This is a simple application of
\ref{equivalent conditions for being fully reducible};  The
irreducible polynomials "$p_i$" in this theorem are in the present
case $p_i=t-b_i$ and,as noticed above, the matrix $V$ appearing in
the statement $(3)$ of \ref{equivalent conditions for being fully
reducible} is exactly the Vandermonde matrix $V_n(b_1,\dots,b_n)$.
The rest is clear.
\end{proof}

Since a $W$-polynomial is of the form $f_\Delta$ for some finite
subset $\Delta \subset K$, the above lemma also shows the strong
relation existing between $W$-polynomials and Vandermonde
matrices.  This leads to the following theorem which shows in
particular, that a companion matrix $C_f$ is
$(S,D)$-diagonalizable if and only if $f$ is a $W$-polynomial.

\begin{theorem}
\label{diagonalization of a companion matrix} Let $f\in R$ be a
monic polynomial of degree $n$. Then the following are equivalent
:
\begin{enumerate}
\item[$i)$] $f$ is a $W$-polynomial.
\item[$ii)$] There exists a P-independent set
$B=\{b_1,b_2,\dots,b_n\} \subset K$ such that $f=f_B$.
\item[$iii)$] There exist $\{b_1,b_2,\dots,b_n\} \subset K$ such
that $V=V_n(b_1,b_2,\dots,b_n)$ is invertible and $$ C_fV =
S(V)\diag(b_1,b_2,\dots,b_n) + D(V)$$
\item[$iv)$] $C_f$ is $(S,D)$-diagonalizable.
\item[$v)$] The left R-module $R/Rf$ is semi-simple with simple
components of dimension 1 over $K$.
\end{enumerate}
\end{theorem}
\begin{proof}
These equivalences are special cases of \ref{equivalent conditions
for being fully reducible} using Lemma \ref{vandermonde and
independence}.
\end{proof}

\begin{remark}
Let us mention that the behaviour here is specific to the left
$R$-module $R/Rf$.  In fact, if $S$ is not onto, even right
modules such as $R/(t-a)R$ need not be semisimple.  Consider for
instance the field $K:=k(x)$ and the $k$-endomorphism $S$ given by
$S(x) = x^2$. If $f(t):= t \in R=K[t;S,D]$ then the $R$-module
$R/fR$ is finitely generated but not artinian (it contains the
descending chain of right $R$-modules $xt^nR + tR$ for $n \in
\mathbb N$) and so cannot be semisimple.
\end{remark}

For the more general case of a matrix $A$ we will assume that the
endomorphism $S$ is an automorphism.   Let us recall that, in this
case, the ring $R=K[t;S,D]$ is in fact a left and right principal
ideal domain.  We will need the following definitions:

\begin{definitions}
For $f,g\in R=K[t;S,D]$ we say that $f$ strongly divides $g$, and
we write $f||g$, if there exists an invariant element $c\in R$
(i.e. $cR=Rc$) such that $f$ left divides $c$ and $c$ left divides
$g$
\end{definitions}

Notice, in particular, that if $f,g\in R$ are such that $f||g$
then $f$ divides $g$ on both sides i.e. $g\in Rf \cap fR$.  In
fact, it is easy to check that the notion of strong divisibility
is left right symmetric.

We can then use the following classical result (Cf. \cite{Co2}).

\begin{lemma}
Let $R$ be a principal ideal domain and let $A$ be an $n \times n$
matrix with coefficients from $R$.  Then there exist invertible $n
\times n$ matrices $P$ and $Q$ such that the matrix

$$PAQ = \diag(e_1,e_2,\dots,e_n)$$

where $e_i$ strongly divides $e_{i+1}$ for $1\le i \le n-1$.
\end{lemma}

A matrix $A\in M_n(K)$ determines a left $R=K[t;S,D]$-module
structure on the space of rows $K^n$. More precisely this
structure is given by $t.\underline v = S(\underline v)A +
D(\underline v)$ (in other words the action of $t$ is given by the
map $T_A$ defined before Definition \ref{diagonalization,
triangulation}). We thus have an exact sequence of left
$R$-modules:

$$0 \longrightarrow R^n \stackrel{tI - A}{\longrightarrow} R^n
\stackrel{\varphi}\longrightarrow K^n \longrightarrow 0$$

where $\varphi$ is the left $R$-morphism sending the unit vectors
of $R^n$ to the unit vectors of $K^n$.  The above lemma shows that
there exist matrices $P,Q \in GL_n(R)$ such that $P(tI - A)Q =
\diag(e_1,e_2,\dots,e_n)$. Remarking that if $e=1$ then $R/eR =
0$, we get after reindexing the $e_i$'s if necessary an
isomorphism of left $R$-modules
\begin{equation}
\label{decomposition of K^n with invariant factors}
_RK^n\cong\bigoplus_{i=1}^r \frac{R}{Re_i}\quad \mathrm{for}\
r\leq n
\end{equation}

The elements $e_i$ in this decomposition are called the invariant
factors. We are now ready for the characterization of an
$(S,D)$-diagonalizable matrix. The last invariant factor "$e_r$"
will play a very important role in the characterization of
$(S,D)$-diagonalizability and triangulability.

\begin{theorem}
\label{diagonalization general case} Let $K,S,D$ be a division
ring, an automorphism and a S-derivation of $K$, respectively. A
matrix $A \in M_n(K)$ is $(S,D)$-diagonalizable if and only if its
last invariant factor is a W-polynomial.
\end{theorem}
\begin {proof}
We continue using the above notations in particular $_RK^n$ is
decomposed as in \ref{decomposition of K^n with invariant
factors}. Since the action of $t.$ is determined by $A$ on $K^n$
and by the $C_{e_i}$ on $R/Re_i$ it then follows from classical
facts (Cf. \cite{L}) that there exists an invertible matrix $P$
such that
\begin{equation}
\label{similarity with invariant factors}
 S(P)AP^{-1} + D(P) = \diag(C_{e_1},C_{e_2}, \dots,C_{e_r})
\end{equation}

It is easy to check that, if the matrices $C_{e_i}$'s are
$(S,D)$-diagonalizable then the matrix $\diag(C_{e_1},C_{e_2},
\dots,C_{e_r})$ is $(S,D)$-diagonalizable.  Conversely: assume
that the matrix $\diag(C_{e_1},C_{e_2}, \dots,C_{e_r})$ is
$(S,D)$-diagonalizable.   This matrix represents the action of
$t.$ (left multiplication by $t$) on $_RK^n\cong \bigoplus_{i=1}^r
\frac{R}{Re_i}$. Hence there exists a $K$-basis
$\{u_1,u_2,\dots,u_n\}$ of $K^n$ consisting of eigenvectors for
the action of $t.$.  We thus have, for
 $l\in \{1,2,\dots,n\}$, $t.u_l=\alpha_lu_l$
 for some $\alpha_l\in K$.
 Decomposing each $u_l$ according to the direct sum $\bigoplus_{i=1}^r
 \frac{R}{Re_i}$,we can write $u_l = \sum_{j=1}^ru_{l,j}$.  It is then easy to
check that for all $j=1,\dots ,r$, the set $\{u_{l,j} \vert \;
i=1,\dots,n\}$ form a generating family of elements of $R/Re_j$
which are eigenvectors for the action of $t$.  We can thus extract
from this family a basis for $R/Re_i$ consisting of eigenvectors.
The union of these families then gives a basis of $K^n$ whose
elements are eigenvectors.  It is now clear that $A$ is
$(S,D)$-diagonalizable if and only if the matrices $C_{e_i}$'s are
$(S,D)$-diagonalizable. Theorem \ref{diagonalization of a
companion matrix} shows that this is the case if and only if the
polynomials $e_1,e_2,\dots,e_r$ are W-polynomials. Since we know
that $e_i$ divides $e_{i+1}$ the conclusion of the theorem follows
from Corollary \ref{factors of W polynomials}.
\end{proof}

The above theorem was obtained using other techniques  by G.
Cauchon in the special case when $S=\id$ and $D=0$ (in particular
Cauchon didn't use the Vandermonde matrices and uses a different
technique of diagonalization).

Let us now come to triangulation.  The expected result holds: a
square matrix $A$ is triangularizable if and only if the last
invariant factor of $A$ is a product of linear factors.  As in the
case of diagonalization we will reduce the problem to the case of
a companion matrix.

\begin{prop}
\label{triangulation of a companion matrix} Let $f\in R=K[t;S,D]$
be a monic polynomial of degree $n$.  The following are eqivalent
:
\begin{enumerate}
\item[i)] $C_f$ is $(S,D)$-triangularizable.
\item[ii)] There exists a chain of left $R$-modules of $R/Rf$

$$0=V_0\lneqq V_1 \lneqq \dots \lneqq V_{n-1} \lneqq V_n = R/Rf .$$
\item[iii)] There exists $g_1,g_2, \dots,g_{n-1} \in R$ such that
$$Rf \subsetneq Rg_1 \subsetneq \dots \subsetneq Rg_{n-1}\subsetneq R .$$
\item[iv)] $f$ is a product of monic linear polynomials.
\end{enumerate}
\end{prop}
\begin{proof}
i) $\longrightarrow$ ii) $C_f$ represents the left multiplication
 $t.: R/Rf \longrightarrow R/Rf$ in the basis $1,t,\dots,t^{n-1}$.
Since $C_f$ is $(S,D)$-triangularizable one can find
$v_1,\dots,v_n$ a $K$-basis of $R/Rf$ such that $t.v_{i} \in Kv_1
+ \cdots + Kv_i$.  In particular, for any $i=1,\dots, n$, the left
$K$-vector space $V_i = Kv_1 +\cdots + Kv_i$ is in fact a left
$R$-module. From this we conclude that these modules satisfy the
required property.

ii) $\longrightarrow$ iii) Thanks to the lemma \ref{companion
matrix of a product} we can find $g_1,\dots,g_n \in R$ such that
$V_i = Rg_i/Rf$. The properties of the $V_i$'s give the required
inclusions between the $Rg_i$'s.

iii) $\longrightarrow$ iv Since $\deg f=n$ and the inclusions are
strict we must have $\deg g_i=n-i$ for $i=1,\dots,n-1$ and we
conclude easily.

iv) $\longrightarrow$ i) Let us write $f(t)= (t-a_1)\dots
(t-a_n)$.  Lemma \ref{companion matrix of a product} (2) shows
that $C_f$ is $(S,D)$-triangularizable.
\end{proof}

We are now ready to present the general case of the criterion for
(upper) triangulation.  For a square matrix $A \in M_n(K)$ we
denote, as in Theorem \ref{diagonalization general case}, by
$e_1,\dots,e_r$ the invariant factors of $A$.  Recall that we have
$e_1 \vert \vert e_2 \vert\vert \cdots \vert \vert e_r$, which
means that there exist invariant polynomials $c_r,\dots,c_1$ such
that $e_i \vert c_i \vert e_{i+1}$.

\begin{theorem}
\label{triangulation general case} Let $K,S,D$ be a division ring
an automorphism and a S-derivation of $K$, respectively. Let $A\in
M_n(K)$ be a square matrix, then $A$ is $(S,D)$-triangularizable
if and only if the last invariant factor $e_r$ is a product of
monic linear polynomials.
\end{theorem}
\begin{proof}
Assume that $e_r$ is a product of linear polynomials.  The fact
that $R$ is a U.F.D. and since we have $e_1 \vert\vert e_2
\vert\vert \cdots \vert\vert e_r$, it is clear that $e_1,\dots ,
e_r $ are also product of linear polynomials.  Proposition
\ref{triangulation of a companion matrix} makes it clear that the
matrices $C_{e_i}$ are all triangularizable. Thanks to equation
\ref{similarity with invariant factors}, we know that $A$ is
similar to $\diag(C_{e_1},\dots,C_{e_r})$ and the result is now
clear. Conversely assume that $A \in M_n(K)$ is triangularizable.
$K^n$ is a left $R$-module via the action
$t.\underline{v}:=S(\underline{v})A + D(\underline{v})$ and let
 $v_1,\dots,v_n$ be a basis of $K^n$ such that, for all $i\in\{1,\dots,n\}$
 $t.v_i = \sum_{j=1}^i\alpha_{ij}v_j$.  Decomposing each $v_i$
 according to the isomorphism \ref{decomposition of K^n with invariant
 factors} we get $v_i =\sum_{k=1}^r v_{ik}$ and so we obtain on one
 hand $t.v_i = \sum_{j=1}^i\alpha_{ij}v_j
 = \sum_{k=1}^r(\sum_{j=1}^i\alpha_{ij}v_{jk})$ and on the other
 hand we have $t.v_i = t.\sum_{k=1}^r v_{ik}=\sum_{k=1}^r
 t.v_{ik}$.  Since $R/Re_k$ is stable by the action of $t.$ and
 the decomposition in \ref{decomposition of K^n with invariant
 factors} is direct we get, for $k\in\{1,\dots,r\}, \;
 t.v_{ik}= \sum_{j=1}^i\alpha_{ij}v_{jk}$.  Let us now observe
 that, for $k=1,\dots,r$, $\{v_{ik} \vert \, i=1,\dots,n\}$ is a generating set
 for $R/Re_k$ as left $K$ vector space.  It is now easy to check that
 one can extract a basis $B_k$ from
 this generating set such that the matrix representing $t.\vert_{R/Re_k}$ in
 the basis $B_k$ is triangular.  Proposition \ref{triangulation of a companion
 matrix} then shows that the $e_k$'s are product of
 linear polynomials.
\end{proof}

%---------------------------------------------------
%---------------------------------------------------
\section{eigenvalues}

In this section we will give some basic facts on eigenvalues of
matrices over division rings.  We will again assume that $S$ is an
automorphism of the division ring $K$.
%The reason for this
%assumption is that we will be able in this case to speak about the
%right and left eigenvalues and define a spectrum
We have seen in the preceding section (see also the paragraph
preceding definition \ref{diagonalization, triangulation}) how to
associate with every matrix $A\in M_{n\times n}(K)$ a structure of
left $R$-module on $K^n$ or equivalently how to define a pseudo
linear transformation $T_A: K^n \longrightarrow K^n$.  Since $S$
is assumed to be an automorphism, the concept defined so far must
be symmetric. The aim of the next lemmq is to examine more closely
this symmetry.
\begin{lemma}
\label{comparison left and right}
\begin{enumerate}
\item $\delta:=-DS^{-1}$ is a right $S^{-1}$-derivation; i.e.
$\delta(ab)=\delta(a)S^{-1}(b)+a\delta(b)$ and $R=K[t;S,D]$ is a
left and right principal ideal domain. The elements of $R$ can be
written in the form $\sum_{i=0}^nt^ia_i$ with the commutation rule
$at=tS^{-1}(a) -DS^{-1}(a)$ for any $a\in K$.
\item We have
$\D^{S,D}(a):=\{a^c:=S(c)ac^{-1} +D(c)c^{-1}\,\vert \,c\in
K\setminus \{0\} \} = \D^{-DS^{-1},S^{-1}}(a):=\{ {^ca}:=
caS^{-1}(c^{-1})+c(-DS^{-1}(c^{-1})) \,\vert \, c \in K \setminus
\{0\} \}$.
\item If $A \in M_n(K)$, we can define a structure of right
$R$-module on the set ${^nK}$ of columns via $u.t:=L_A(u):=
AS^{-1}(u)-DS^{-1}(u)$ where $u \in {^nK}$.
\item If $A\in M_n(K)$ the left
$R$-module $K^n$ and the right $R$-module ${^nK}$ induced by $A$
gives rise to the same invariant factors (up to similarity).  i.e.
$K^n \cong \bigoplus_{i=1}^r R/Re_i  \Leftrightarrow {^nK}\cong
\bigoplus_{i=1}^r R/e_iR$.
\end{enumerate}
\end{lemma}
\begin{proof}
(1) This is standard and easy to prove.

(2) It suffices to check that for $c\in K \setminus \{0\}$ we have
$ {^ca}=a^d$ where $d=S^{-1}(c)$.

(3) Let us compute, for  $\alpha \in K$ and $u\in {^nK}, \,
L_A(u\alpha)=
AS^{-1}(u\alpha)-DS^{-1}(u\alpha)=AS^{-1}(u)S^{-1}(\alpha)-D(S^{-1}(u)S^{-1}(\alpha))
=
AS^{-1}(u)S^{-1}(\alpha)-uDS^{-1}(\alpha)-D(S^{-1}(u))S^{-1}(\alpha)=
L_A(u)S^{-1}(\alpha)+u(-DS^{-1})(\alpha)$.  This shows that
$(u\alpha).t =
(u.t)S^{-1}(\alpha)+u(-DS^{-1})(\alpha)=u.(tS^{-1}(\alpha)-(DS^{-1})(\alpha))
=u.(\alpha t)$.  The rest is clear.

(4) This is due to the fact that the invariant factors are
obtained from $tI - A \in M_n(R) $ using elementary
transformations on rows and columns and hence depend only on $A$.
\end{proof}

\begin{definition}
For $A \in M_{n\times n}(K), \, \alpha,\beta \in K, \; v\in
K^n\setminus \{(0,\dots,0)\}$ and $u\in\;^nK\setminus
\{(0,\dots,0)^t\} $, we say that:
\begin{enumerate}
\item $\alpha$ is a left eigenvalue of $A$ associated to $v$ if
$$T_A(v)=\alpha v$$.
\item $\beta$ is a right eigenvalue of $A$ associated to $u$ if
$$L_A(u)=u\beta$$
\end{enumerate}
\end{definition}

We will denote $\lspec(A)$ and $\rspec(A)$ the sets of left and
right eigenvalues of a matrix $A$; $\Spec(A)$ will denote the
union of left and right eigenvalues.

In the next proposition we collect a few elementary properties of
the left and right eigenvalues.

\begin{prop}
Let $A$ be a matrix in $M_n(K)$.  Then,
\begin{enumerate}
\item $\lspec(A), \, \rspec(A),\, \Spec(A)$ are closed under
$(S,D)$-conjugation.
\item If $P \in GL_n(K)$,
$$\lspec(A)=\lspec(A^P), \, \rspec(A)=\rspec(A^P), \,
\Spec(A)=\Spec(A^P)\,.$$
\item Left eigenvectors corresponding to non $(S,D)$-conjugate
left eigenvalues are left linearly independent.
\item Right eigenvectors corresponding to non $(S,D)$-conjugate
right eigenvalues are right linearly independent.
\item If $\alpha \in \lspec(A)$ and $\beta \in \rspec(A)$ are not $(S,D)$-conjugate and
$v=(v_1,\dots ,v_n)\in K^n,\,u=(u_1,\dots,u_n)^t \in {^nK} $ are
the associated eigenvectors then $v.u:=\sum_{i=1}^n v_iu_i=0$.
\end{enumerate}
\end{prop}
\begin{proof}
(1) Assume $\alpha \in \lspec(A)$ and let $v\in K^n$ be an
eigenvector for $\alpha$.  We thus have $T_A(v)=\alpha v$.  If
$\beta \in K\setminus \{0\}$ we have $T_A(\beta
v)=S(\beta)T_A(v)+D(\beta)v=(S(\beta)\alpha +
D(\beta))v=(\alpha^{\beta})\beta v $.  This shows that
$\alpha^{\beta}$ is also a left eigenvalue and proves that
$\lspec(A)$ is closed under $(S,D)$-conjugation. Similarly, if
$\lambda \in \rspec(A)$, $u\in {^nK}$ and $\gamma \in K \setminus
\{0\}$ are such that $L_A(u)=u\lambda$, one can check that
$L_A(uS(\gamma^{-1}))=uS(\gamma^{-1}) \lambda^\gamma$.

(2) It is easy to verify that for $v\in K^n$ we have
$T_{A^P}(v)P=T_A(vP)$.  From this one deduces that if $\lambda \in
K$ is such that $T_{A^P}(v)=\lambda v$ then $T_A(vP)=\lambda vP$;
This shows that $\lspec(A^P)\subseteq \lspec(A)$.  The reverse
inclusion follows since $P\in GL_n(K)$. Similar computations lead
to $\rspec(A)=\rspec(A^P)$

(3),(4) and (5) are easy to prove and can be found in \cite{L},
Proposition $4.13$.
\end{proof}

As in the case when $K$ is a commutative field and $S=id.,D=0$ we
will now show that the eigenvalues are exactly the roots of some
monic polynomials.  In the classical case the last invariant
factor is the minimal polynomial.  This polynomial is unique.  In
our case the last invariant factor is only defined up to
similarity. In Lemma \ref{roots of similar polynomials} we will
compare the roots of similar polynomials. First let us recall that
$f,g\in R $ are said to be similar, denoted $f\sim g$, iff $R/Rf
\cong R/Rg$ if and only if $R/fR \cong R/gR$.    For a polynomial
$f\in R=K[t;S,D]$, we continue to denote $V(f)$ the set of its
right roots i.e. $V(f)=\{a\in K \, \vert \, f\in R(t-a) \}$.
Similarly we will denote $V'(f)$ the set of left roots of $f$ i.e.
$V'(f)=\{a\in K \,|\, f\in (t-a)R \}$.

\begin{lemma}
\label{roots of similar polynomials}
 Let $f,g$ be similar elements in $R$.  Assume that $R/Rf
\stackrel{\gamma} \cong R/Rg: 1+Rf \mapsto q+Rg$ then
$V(f)=\phi_q(V(g))$.
\end{lemma}
\begin{proof}
Since $\gamma$ is well defined, there exists $q'\in R$ such that
$fq=q'g$.  The map $\gamma$ being onto, we must have $Rq+Rg=R$. In
particular, $V(q)\cap V(g)=\emptyset$.  So if $x \in V(g)$, we
have $x\in V(fq)\setminus V(q)$ and the formula \ref{eq1} implies
that $\phi_q(x)\in V(f)$. We thus conclude that
$\phi_q(V(g))\subseteq V(f)$.  Similarly if
$\gamma^{-1}(1+Rg)=p+Rf$, we must have $\phi_p(V(f))\subseteq
V(g)$.  We also have $qp \in 1+Rf$ and this implies that
$\phi_{qp}$ is the identity on $V(f)$.   It is also easy to check
that $\phi_{qp}=\phi_q \circ \phi_p$ (Cf. \cite{LL5}).  We thus
get:
$$
V(f)=\phi_{qp}(V(f))=\phi_q(\phi_p(V(f)))\subseteq
\phi_q(V(g))\subset V(f) \, .
$$
This yields the result.
\end{proof}

\begin{corollary}
\label{roots of similar polynomials are conjugete}
 If $f,g\in
R=K[t;S,D]$ are similar there exist $p,q\in R$ such that $V(g)\cap
V(q)=V(f)\cap V(p)=\emptyset$ and $V(f)=\{\alpha^{q(\alpha)}\,|\,
\alpha \in V(g)\}$ and $V(g)=\{\beta^{p(\beta)}\,|\,\beta \in
V(f)\}$.
\end{corollary}

Of course, there exist similar statements for the left roots using
the left analogue of the map $\phi$.

We can now give the analogue of the classical fact that the roots
of the minimal polynomial are exactly the eigenvalues of the
matrix.
\begin{prop}
Let $A\in M_n(K)$ and $\{e_1,\dots,e_r\}$ be a matrix and a
complete set of invariant factors for $A$.  Denote by
$\Delta(e_r)$ the set $\{f\in R \, \vert \, f \sim e_r \} $, then
the following are equivalent:
\begin{enumerate}
\item[i)] $\beta \in \rspec(A)$.
\item[ii)] There exists $\gamma \in K\setminus \{0\}$ such that $\beta^\gamma \in
V(e_r)$.
\item[iii)] There exists a polynomial $e'_r \in \Delta (e_r)$ such that $\beta \in
V(e'_r)$.
\end{enumerate}
Similar statements hold for elements of $\lspec(A)$ and $V'(e_r)$.
\end{prop}
\begin{proof}
$(i) \Rightarrow (ii)$ Assume $u \in\, ^nK \setminus \{0\}$ is
such that $L_A(u)=u\beta $.   This also means that while
considering $^nK$ as a right $R$-module, $u.(t-\beta)=0$. Writing
$u=(u_1+e_1R,\dots,u_r+e_rR)$  according to the decomposition
obtained in Lemma \ref{comparison left and right}, we get that
there exists $i \in \{1,\dots,r\}$ such that $u_i\notin e_iR \ne
0$ but $u_i(t-\beta)\in e_iR$.  We may assume that $deg(u_i)<deg
(e_i)$ and, comparing degrees, we conclude that there exists an
element $\gamma \in K\setminus \{0\}$ such that $u_i(t-\beta)=
e_i\gamma$.  This leads to
$u_iS(\gamma^{-1})(t-\beta^\gamma)=e_i$. Since $e_i$ divides $e_r$
on the right, we do get that $\beta^\gamma\in V(e_r)$.

$(ii) \Rightarrow (iii)$  By hypothesis there exists $\gamma\in
K\setminus\{0\}$ and $g\in R$ such that $g(t-\beta^\gamma)=e_r$.
Right multiplying by $\gamma$ we get
$g(t-\beta^\gamma)\gamma=e_r\gamma$ i.e.
$gS(\gamma)(t-\beta)=e_r\gamma$.  This yields the result since
$e'_r:=e_r\gamma $ is obviously similar to $e_r$.

$(iii) \Rightarrow (ii)$ This is clear from Corollary \ref{roots
of similar polynomials are conjugete}.

$(ii) \Rightarrow (i)$ Since $\beta^{\gamma} \in V(e_r)$, we
easily get that $\beta^{\gamma} \in rspec(A)$ and the fact that
$\rspec(A)$ is closed by $(S,D)$ conjugation implies that $\beta
\in \rspec(A)$.

The statements for $\lspec(A)$ and $V'(e_r)$ are similar using
$T_A$ instead of $L_A$ as well as Lemma \ref{comparison left and
right}.

\end{proof}

We can now conclude:

\begin{corollary}
Let $A$ be a matrix in $ M_n(K)$ and $\{e_1,\dots,e_r\}$ be a
complete set of invariant factors for $A$ such that $e_1 \vert
\vert e_2 \dots \vert\vert e_r$. Then
\begin{enumerate}
\item $$ \lspec (A) = \cup_{f\in \Delta(e_r)} V'(f) \, .$$
\item $$ \rspec (A) = \cup_{f\in \Delta(e_r)} V(f) \, .$$
In particular, if $\Gamma_{r}:=\{q\in R \,\vert \,Rq+Re_r=R
\;{and}\;\deg q < \deg e_r\}$ then $rspec (A) = \bigcup_{q\in
\Gamma_r} \phi_q(V(e_r))$.
\end{enumerate}
\end{corollary}

\begin{corollary}
Let $A$ be a matrix in $M_n(K)$.  The number of non
$(S,D)$-conjugate elements in $\Spec (A)$ is bounded by $deg
(e_r)$.
\end{corollary}
\begin{proof}
Notice that if $f\in \Delta (e_r)$, Corollary \ref{roots of
similar polynomials are conjugete} shows that the conjugacy
classes intersecting $V(f)$ also intersects $V(e_r)$.  Hence the
($S,D$) conjugacy class intersecting $\rspec (A)$ also intersects
$V(e_r)$.  Similarly the ($S,D$) conjugacy classes intersecting
$\lspec (A)$ also intersects $V'(e_r)$.  Now, Corollary
\ref{Conjugacy classes intersecting roots} shows that the number
of ($S,D$)-conjugacy classes intersecting $\Spec (A)$ is bounded
by $deg (e_r)$.
\end{proof}

\section{G-algebraic sets and G-polynomials}

In this section we will restrict our attention to the case when
$S=id.$ and $D=0$.  $K$ will stand for a division ring, $G$ will
denote a group of automorphisms of $K$ and $K^G:= \{x \in K \vert
\sigma(x) = x \; \forall \sigma \in G \}$.

\begin{definition}
\label{G-algebraic} A subset $\Delta \subseteq K$ is $G$-algebraic
if there exists a monic polynomial $f\in K^G[t]$  such that
$f(x)=0$ for all $x\in \Delta$.  The monic polynomial in $ K^G[t]
$ of minimal degree annihilating $\Delta$ is denoted $f_{\Delta
,G}$.  Polynomials of the form $f_{\Delta ,G}$ will be called
$G$-polynomials.  In particular, if $G=\{Id.\}$ we find back the
notion of an algebraic set in the sense defined in Wed1
(\cite{LL5}).
\end{definition}

It will sometimes be useful to denote the unique monic least left
common multiple of a set $\Gamma$ of (monic) polynomials by
$\Gamma_\ell$.  Of course every $G$-algebraic set is algebraic;
the next proposition gives characterizations of $G$-algebraic
sets.

\begin{prop}
\label{characterization of G-algebraic sets} With the above
notations, the following are equivalent:
\begin{enumerate}
\item[i)] $\Delta$ is $G$-algebraic.
\item[ii)] $\bigcup_{\sigma \in G}\;\sigma(\Delta)$ is algebraic.
\item[iii)] $\Delta$ is algebraic and for all $a \in\Delta$,
$\{\sigma(a) \vert \sigma \in G\}$ is algebraic.
\item[iv)] $\Delta$ is algebraic and if
$\{a_1,a_2,\dots,a_n\}$ is a $P$-basis for $\Delta$ then $\{a_i\}$
is $G$-algebraic for $ 1\le i \le n$.
\item[v)] There exists a left common multiple of the set
$\{t - \sigma(a)\, \vert \,\sigma \in G , \; a \in \Delta\}$
\end{enumerate}
\end{prop}
\begin{proof}
i) $\Longrightarrow$ ii) If $f\in K^G[t]$ is such that $f(\Delta)
=0$ then $f(\Delta^\sigma) = 0$ for all $\sigma \in G$.  Hence
$f(\cup_{\sigma \in G}\;\sigma(\Delta)) = 0$.

\noindent ii) $\Longrightarrow$ iii) Since $\Delta \subseteq
\cup_{\sigma \in G}\; \sigma(\Delta)$, we have that $\Delta $ is
algebraic.  Similarly for all $a\in \Delta$, $G.a:= \{\sigma(a)
\vert \sigma \in G\}\subseteq \cup_{\sigma \in G}\;\sigma
(\Delta)$, hence $G.a$ is algebraic and its minimal polynomial is
precisely the monic generator of the left ideal $\bigcap_{\sigma
\in G}R(t-\sigma(a)) \ne 0$.  In other words,
$f_{G.a}=\{t-\sigma(a)\vert \sigma \in G\}_\ell \in K^G[t]$.

\noindent iii) $\Longrightarrow$ iv) This is obvious.

\noindent iv) $\Longrightarrow$ v) Let $\{a_1,a_2,\cdots,a_n \}$
be a $P$-basis for for $\Delta$ and define $f_i$ to be the left
common multiple of the set $ \{t-\sigma(a_i)\vert \sigma \in G
\}$. Then $f_i^\sigma = f_i$ , i.e. $f_i \in K^G[t]$ for all $i\in
\{1,2,\cdots,n\}$. Hence we have $f:=\{f_i \vert i=1,2, \cdots,n
\}_\ell = \{t-\sigma(a)\vert \sigma \in G \;,\; a \in
\{a_1,a_2,\cdots,a_n\}\}_\ell \in K^G[t]$. But $a \in \Delta$
implies that $t-a$ divides on the right $\{t-a_i\vert i \in
\{1,2,\cdots,n\}\}_\ell$ which itself divides $f$ on the right.
Since $f \in K^G[t]$ we thus get that $f$ is a left common
multiple of the set $ \{t-\sigma(a)\vert \sigma \in G \; , \; a
\in \Delta \}$.

\noindent iv) $\Longrightarrow$ i)  This is left to the reader.
\end{proof}

\begin{remarks}
\begin{enumerate}
\item[a)] Of course if $G$ is a finite group then every algebraic
set is $G$-algebraic.
%*******
%compare the rank of an algebraic set and the rank of its $G$-closure
%*******
\item[b)] Notice that in the case when $K$ is commutative, a
$G$-algebraic set must be finite.
\item[c)] Part $iv)$ of the above proposition explains why we will be
mainly concerned with $G$-algebraic sets of the form
$\{\sigma(a)\vert \,\sigma \in G \}$ for some $a\in K$; this set
will be denoted by $G.a$.

\item[d)] If $\Delta $ is an algebraic set and $\sigma$ is an
automorphism then $\sigma(\Delta)$ is also algebraic its minimal
polynomial is $\sigma (f_\Delta)$ where we assume that $\sigma$
has been extended to $K[t]$ by putting $\sigma (t)=t $.  In
particular we get that $\rk\Delta = \rk \sigma(\Delta)$.
\end{enumerate}
\end{remarks}

\begin{corollary}
\label{splitting of G-polynomials} Any $G$-polynomial
$f=f_{\Delta,G}$ factorizes linearly: $f=(t-b_1)\cdots (t-b_n)$ in
$K[t]$. Moreover any root of $f$ is conjugated to some $b_i$'s and
these $b_i$'s are conjugated to elements in $\bigcup_{\sigma \in
G} \sigma(\Delta)$.
\end{corollary}
\begin{proof}
These are obvious consequences of the above  proposition and of
our earlier results in \cite{LL5}.
\end{proof}

\begin{examples}
{\rm
\begin{enumerate}
\item[a)] Let $G$ be the set of all inner automorphisms of $K$ i.e.
$G=\{I_x \vert\, x\in K^* \}$.  Then $K^G = Z(K)$ the center of
$K$. An element is then $G$-algebraic if it is algebraic over the
center $Z(K)$.  In particular the above corollary gives back the
Wedderburn classical theorem: If an element $a$ of a division ring
$K$ is algebraic over the center $Z(K)$ then its minimal
polynomial factorizes in $K[t]$ into linear factors of the form
$t-b$ where $b\in K$ is conjugate to $a$.
%********
%can we say something about cyclic permutations in our general setting ?
%********
\item[b)] Let $D$ be a division subring of $K$ and put $L=C_K(D)$
the centralizer of $D$ in $K$.  Then $L=K^G$ for $G=\{I_x \vert x
\in D^*\}$ hence an element $a\in K$ is algebraic over $L$ if and
only if it is $G$-algebraic.  In this case, the above corollary
shows that its minimal polynomial over $L$ factorizes linearly in
$K[t]$.  Notice that in the case when $K$ is finitedimensional
over its center $Z(K)$ then every subdivision ring $L$ such that
$Z(K) \subseteq L \subseteq K$ is such that $L=C_K(C_K(L))$ and
the conclusion applies.
\item[c)] If $K$ is commutative and $G$ is a subgroup of
automorphisms of $K$, an element $a \in K$ is algebraic over
$L=K^G$ if and only if the set $\{\sigma(a)\vert \sigma \in G\}$
is finite.  We also get back the classical fact on galois
extensions: every such extension is normal.
\end{enumerate}
}
\end{examples}
\begin{theorem}
\label{G1} Let $G$ be a group of automorphisms of $K$, and suppose
that $a\in K$ is algebraic over $K^G$.  Define $G_a:=\{\sigma \in
G \,\vert \, \sigma(a) \in \Delta (a)\}$, where $\Delta (a)=\{a^x
\, \vert \, x\in K \setminus \{0\} \} $. Then:
\begin{enumerate}
\item[a)] $G_a$ is a subgroup of $G$.

\item[b)] For any $\sigma \, , \, \tau \in G$ we have
$\sigma G_a = \tau G_a$ (resp. $G_a\sigma =G_a\tau $) if and only
if $\Delta (\sigma(a))= \Delta (\tau (a))$ (resp. $\Delta
(\sigma^{-1}(a))=\Delta (\tau^{-1}(a))$).
\item[c)] $G_a$ is of finite index in $G$.

\item[d)] The decomposition of $G$ into its right cosets modulo
$G_a$ corresponds to the decomposition of $G.a$ into conjugacy
classes.  More precisely if $G = \bigcup_{i=1}^{n}\sigma_iG_a$ is
the decomposition of $G$ into its right cosets modulo $G_a$ then
$G.a = \bigcup_{i=1}^{n}\sigma_i(G_a.a)$ is the decomposition of
$G.a$ into conjugacy classes.

\item[e)] $\rk(G.a) = \deg f_{a,G}=(G:G_a)\rk G_a.a=(G:G_a)\deg f_{a,G_a}
= (G:G_a)\dim_CYC$ where $Y\subseteq K\setminus \{0\}$ is such
that $G_a.a=a^Y$. More precisely, if $\{y_1,y_2,\cdots,y_n\}$ is a
maximal $C$-independent set in $Y$ then $\sigma (a^{y_j})$ is a
$P$-basis for $G.a$.

\item[f)] If $G_a=\{Id.\}$ then $G_{int.}:=
\{\sigma \in G \,|\, \sigma \,{\rm is \; inner} \}=\{Id.\}$.
Moreover, if $\sigma$ and $\tau$ are different elements in $G$,
then $\sigma (a)$ and $\tau (a)$ belong to different conjugacy
classes and $G_a$ is full.
\end{enumerate}
\begin{proof}
a) This is left to the reader.

\noindent b) Suppose $\sigma G_a=\tau G_a$.  We can write $\sigma
= \tau g_1$ for some $g_1 \in G_a$.  The definition of $G_a$ shows
that there exists $x_1 \in K $ such that $g_1(a)=a^{x_1}$.  For
$y\in K$ we then have $\sigma(a)^y=\tau(g_1(a))^y=\tau(a^{x_1})^y
=(\tau (a)^{\tau(x_1)})^y=\tau (a)^{y\tau(x_1)}$.  This shows that
$\Delta (\sigma (a))\subseteq \Delta(\tau(a))$.  The reverse
inclusion is proved similarly.

\noindent The proof of sufficiency of the condition as well as the
proof of the analogue left-right statements are left to the
reader.

\noindent c) Since $G.a$ is algebraic it can only intersects a
finite number of conjugacy classes i.e. the number of conjugacy
classes of the form $\Delta (\sigma(a) )$ where $\sigma \in G$ is
finite.  Part b) above enables us to conclude.

\noindent d)  This is easily deduced from $b)$ above.

\noindent e) This is a direct consequence of $d)$ above using
results from \cite{LL2}.

\noindent f) Theses are easy consequences the definitions.
\end{proof}
\end{theorem}

Let us remark that the subgroup $G_a$ contains the subgroup
$G_{int}$ of all the inner automorphisms.

\begin{example}
{\rm The condition $(G:G_a)<\infty$ is not sufficient for $a$  to
be $G$-algebraic: for instance if $G=G_{int}$, then $K^G=Z(K)$,
the center of $K$ and $G=G_a$ for any $a\in K$ but of course $a$
 is not necessarily algebraic over $Z(K)$.} \end{example}

Before giving necessary and sufficient conditions for $a$ to be
$G$-algebraic let us recall that a subset of a conjugacy class
$\Delta(a)$, say $a^Y$, is algebraic if and only if the right
$C(a)$-vector space $YC(a)$ generated by $Y$ over the centralizer
of $a$ is finitedimensional. (Cf. Prposition 4.2 in \cite{LL2})
\begin{prop}
\label{NSC for G.A to be algebraic}  Let $a$ be an element of $K$
and $Y$ a subset of $K\setminus\{0\}$ such that
$G_a.a=\{a^y\,|\,y\in Y\}$.  Then $a$ is $G$-algebraic if and only
if the right $C(a)$-vector space generated by $Y$ is
finitedimensionnal and $(G:G_a)<\infty$.
\end{prop}
\begin{proof}
If $G.a$ is algebraic we have seen in Theorem \ref{G1} that
$(G:G_a)<\infty$.  On the other hand since $G_a.a\subseteq G.a$,
it is clear that $G_a.a$ is an algebraic subset contained in
$\Delta (a)$.  This implies that the $C(a)$-right vector space
generated by $Y$ is finitedimensional.

\noindent Conversely, Suppose that $(G:G_a)<\infty$ and let
$\sigma_1,\dots,\sigma_l$ be such that $G=\cup_{i=1}^l\sigma _iG_a
$, then
$G.a=\cup_{i=1}^l\sigma_iG_a.a=\cup_{i=1}^l\sigma_i(a)^{\sigma_i(Y)}$
is the decomposition of $G.a$ into conjugacy classes.  It is easy
to check that, for any $i=1,\dots,l$,
$dim_{C(a)}YC(a)=dim_{C(\sigma_i(a))}(\sigma_i(Y)C(\sigma_i(a) )$.
Since $dim_{C(a)}YC(a)<\infty$, we conclude that the subsets
$\sigma_iG_a.a$ are algebraic for $i=1,\dots,l$.  From this and
the decomposition of $G.a$ given above we get the result.
\end{proof}

We will end this section with some results about the
irreducibility of a $G$-polynomial.  First let us notice that a
$G$-polynomial is not always irreducible:

\begin{example}
\label{reducible minimal polynomial of a G algebraic set}

{\rm Let $K= \mathbb H$, the real quaternions and
$G=\{id.,Int(i)\}$, then $K^G= \mathbb C$. Consider $a=j ,\; G.a
=\{j,j^{i}\}$ is algebraic with minimal polynomial $t^2+1\in
\mathbb C[t]$.  Since $t^2+1=(t+i)(t-i)$ we conclude that the
$G$-polynomial $t^2+1$ is reducible in $K^G[t]$. }
\end{example}

Let us recall, from our earlier work, the following definition:
\begin{definition}
An algebraic set $\Delta \subseteq K$ is said to be full if
$V(f_\Delta)=\Delta$.
\end{definition}

\begin{prop}
\label{irreducibility 1} Let $a\in K$  be a $G$-algebraic element
such that $\Delta:=G.a$ is full then $f_\Delta$ is irreducible in
$K^G[t]$.
\end{prop}
\begin{proof}
Assume $f_\Delta = gh$ in $K^G[t]$.  If $\deg h>0$ then, since
$f_\Delta$ is a $W$-polynomial, we get that $V(h) \ne \emptyset$.
Now if $x\in V(h)$, then $x\in V(f_\Delta)=\Delta$, where the last
equality comes from the hypothesis that $G.a$ is full. Since $h\in
K^G[t]$ we have, for any $\sigma \in G ,\;
0=\sigma(h(x))=h(\sigma(x))$.  We thus get that $h(G.x)=0$. Now
writing $x=\tau(a)$ for some $\tau$ in $G$, we easily get that
$G.x=G.a=\Delta$ and hence, $h(\Delta)=0$. This shows that
$h=f_\Delta$.
\end{proof}

\begin{remark}
\label{irreducibility doesn't imply fullness} The above sufficient
condition for irreducibility in $K^G[t]$ of a minimal polynomial
of a $G$-algebraic set is not necessary, i.e. a $G$-algebraic set
$\Delta$ such that $f_\Delta$ is irreducible in $K^G[t]$ is not
necessarily full. Indeed, consider $K=\mathbb H_{\mathbb Q}$ the
quaternions over the rational numbers,
$G=\{Id.,Int(i)\},\;K^G=\mathbb Q(i)$ and $a=i+j$. Then
$G.a=\{i+j,i-j\}$ is algebraic.  $f_{G.a}\in \mathbb Q(i)[t]$ has
degree $2$ and $V(f_{G.a})=\{(i+j)^{\lambda + i \mu}\, \vert
\,\lambda,\mu \in C_{\mathbb H}(i+j)\}$.  This shows that $G.a$ is
not full.  Now, if $f_{G.a}$ has a root in $\mathbb Q(i)$ then
there exists $x \in \mathbb H_{\mathbb Q}$ such that $(i+j)^x\in
\mathbb Q(i)$.  Let us write $(i+j)^x=\alpha + i\beta$ with
$\alpha,\beta \in \mathbb Q$.  Taking traces on both sides of this
equation, we get $\alpha=0$ and looking at norms we then conclude
that $\beta^2=2$.  Since this last relation is impossible we can
conclude that $f_{G.a}$ is irreducible in $\mathbb Q(i)$.
% to prove the last statement remark that 1) conjugation and standard involution
% commute.  This shows that the norm and the trace of a quaternion are preserved
% by conjugations.  So if $\Delta (i+j) \cap \mathbb Q(i)$ is not empty
% say $(i+j)^x=\alha + i\beta$ where $x\in \mathbb H$ and $\alpha,\beta \in \mathbb C$
% then $0=2\alpha$ and $2=(\beta)^2$.
\end{remark}

The above proposition and theorem \ref{G1} immediately leads to
the following

\begin{corollary}
Assume the group $G_a$ is trivial: $G_a = \{1\}$ then $\Delta =
G.a$ is full and $f_\Delta$ is irreducible in $K^G[t]$.
\end{corollary}

In the same spirit, let us mention the following necessary and
sufficient condition for irreducibility of the minimal
$G$-polynomial associated to a $G$-algebraic set:

\begin{prop}
\label{irreducibility 2} Let $a\in K$ and $\Delta = G.a$ be
algebraic.  Then $f_\Delta$ is irreducible in $K^G[t]$ if and only
if for any $b \in K$ such that $f_\Delta(b)= 0$ we have $f_\Delta
= f_{G.b}$.
\end{prop}
\begin{proof}
Assume $f_\Delta(b) = 0$ then $f_\Delta (G.b) =0$ hence $f_{G.b}$
divides on the right $f_\Delta$ in $K^G[t]$ and the irreducibility
of $f_\Delta$ implies that $f_{G.b}=f_\Delta$.

\noindent Conversely, assume $f_\Delta = gh$ in $K^G[t]$ with $h$
monic and $\deg h \ge 1$, then there exists $x \in \Delta = G.a$
such that $h(x)=0$ and so $h(\Delta)=0$ which shows that
$h=f_\Delta$.

\end{proof}

\author{T. Y. Lam}
\address{Department of Mathematics, University of California,
Berkeley, CA 94720}
%\curraddr{}
\email{lam@math.berkeley.edu}
\thanks{The research reported in this paper was partially supported by a
grant from NSA (T.Y.L.) and a faculty grant from the Universit\'e
d'Artois at Lens (A.L.)}

\author{A. Leroy}
\address{Department of Mathematics, Universit\'e d'Artois, 62307 Lens =
Cedex, France}
% \curraddr{Department of Mathematics, Universit\'e d'Artois, }
\email{leroy@euler.univ-artois.fr}
% \thanks{The author was supported in part by xxxx.}

%    Information for third author
\author{A. Ozturk}
\address{Institut de Math\'{e}matique, Universit\'{e} de Mons-Hainaut,
B-7000 Mons, Belgique.}
% \curraddr{Institute of Mathematics, Universit\'e de Mons-Hainaut, }
\email{ozturk@umh.ac.be}
% \thanks{The author was supported in part by xxxx.}
%
\end{document}